\documentstyle[11pt]{article}
\def\R{{\rm I}\!{\rm  R}}

\def\eps{\varepsilon}

\newcommand{\nl}{\\ \displaystyle} 
\catcode`@=11
\@addtoreset{equation}{section}
\catcode`@=12
 \newtheorem{thm}{Theorem}[section]

\newtheorem{lem}{Lemma}[section]
\newtheorem{cor}{Corollary}[section]
\newtheorem{defi}{Definition}[section]
\newtheorem{esempio}{Example}[section]

\newtheorem{rem}{Remark}[section]

\newcommand{\dis}{\displaystyle}
\title{Barrier functions for Pucci--Heisenberg operators and applications \thanks{ This work was partially
supported by the PRIN MIUR:Metodi di viscosit\'a, metrici e di teoria del controllo in equazioni alle derivate parziali nonlineari.}}
\author{\\ Alessandra Cutr\`\i\\ {\small Dipartimento di Matematica, 
Universit\`a di Roma, \lq\lq Tor Vergata" 00133 Roma, Italy }\\ 
\\ Nicoletta Tchou\\ {\small IRMAR, Universit\'e de Rennes 1, 35042 Rennes, France }\\ 
} 
\date{}

\begin{document}
\maketitle
\begin{abstract}
The aim of this article is the explicit construction of some barrier
functions ("fundamental solutions") for the Pucci-Heisenberg operators.
Using these functions we obtain the continuity
property, up to the boundary, for the viscosity solution of  fully
non-linear Dirichlet problems on the Heisenberg group,
if the boundary of the domain satisfies some regularity geometrical
assumptions (e.g. an exterior Heisenberg-ball condition at the
characteristic points).
We point out that the knowledge of the fundamental solutions allows also to
obtain qualitative properties of Hadamard, Liouville and
Harnack type.
\end{abstract}
\section {Introduction}

In this paper we study viscosity solutions to some degenerate elliptic fully nonlinear second order equations modelled on the Heisenberg vector fields. Precisely, we consider equations of the following type:
\begin{equation}\label{0}
F(\xi, D^2_{H^n}u)+H(\xi,\nabla_{H^n}u)=0,
\end{equation}
where the second order term is obtained by a composition of  a fully non linear  uniformly elliptic operator $F$ with ellipticity constants $0<\lambda\leq \Lambda$, such that $F(\xi,0)=0$, with  the degenerate Heisenberg Hessian matrix $D^2_{H^n}u =(X_iX_ju)_{sym}$, where
\begin{equation}\label{heisf}
X_i=\frac{\partial}{\partial \xi_i} + 
2\xi_{i+n}\frac{\partial}{\partial\xi_{2n+1}}\;,\;\;\;
X_{i+n}=\frac{\partial}{\partial \xi_{i+n}} - 
2\xi_i\frac{\partial}{\partial\xi_{2n+1}}
\end{equation}
for  $i=1,\dots ,n$, $\xi$ denotes the generic point of $\R^{2n+1}$ and  $sym$ denotes the symmetrized matrix.\\
 The first order term $H$ depends on the Heisenberg gradient $\nabla_{H^n}u=(X_1u,\dots ,X_{2n}u)=\sigma(\xi)\nabla u$ ($\sigma$ being the matrix whose rows are the coefficients of the vector fields $X_i$) and satisfies, for some suitable constants $K,M>0$, and for some modulus $\omega$:
 \begin{equation}\label{inH}
 \begin{array}{l}
|H(\xi,\sigma(\xi)p)|\leq K|\sigma(\xi)p|+M\,,\\
|H(\xi,\sigma(\xi)p)-H(\eta,\sigma(\eta)p)|\leq \omega(|\xi-\eta|)(1+|p|)\,.
\end{array}
\end{equation}
The hypotheses on $F$ and $H$ imply that  the following  inequality  holds true in the viscosity sense, see e.g. \cite{caff}:
\begin{equation}\label{sigmastimafm-}
\begin{array}{ll}
 \tilde{{\cal  M}}^-_{\lambda,\Lambda}(D^2u)-K|\nabla_{H^n}u|-&M \leq F(\xi,D^2_{H^n}u)+H(\xi,\nabla_{H^n}u)\leq\\
 &\leq  \tilde{{\cal  M}}^+_{\lambda,\Lambda}(D^2u)+K|\nabla_{H^n}u|+M,
 \end{array}
\end{equation}
where the {\it Pucci--Heisenberg operators} $\tilde{\cal  M}^\mp$ are defined by the composition of the Pucci operators ${\cal  P}^\mp$ (see \cite{puc, caff}) with the Heisenberg hessian $D^2_{H^n}$, that is:
\begin{equation}\label{puc}
\begin{array}{l}
\displaystyle
 \tilde{{\cal  M}}^{-}_{\lambda,\Lambda}(D^2u)\dot={\cal  P}^{-}_{\lambda,\Lambda}(D^2_{H^n}u)=- \Lambda\sum_{e_i>0}e_i-\lambda \sum_{e_i<0}e_i,\\
 \displaystyle
  \tilde{{\cal  M}}^{+}_{\lambda,\Lambda}(D^2u)\dot={\cal  P}^{+}_{\lambda,\Lambda}(D^2_{H^n}u)=- \Lambda\sum_{e_i<0}e_i-\lambda \sum_{e_i>0}e_i,
  \end{array}
 \end{equation}
 where $e_i$ denote the eigenvalues of $D^2_{H^n}u$.

For $\lambda=\Lambda$ all the operators $F$ coincide with  the Heisenberg laplacian $-\Delta_{H^n}\cdot \dot= -\mbox{tr}(D^2_{H^n}\cdot)$. In the linear case, the  operator $F=L_A=-\mbox{tr}(AD^2_{H^n}\cdot)$, where $A=P^{T}P$ is a symmetric $2n\times 2n$ matrix with eigenvalues in between $[\lambda\,,\,\Lambda]$, arises in the 
stochastic theory as the infinitesimal generator of the degenerate 
 diffusion process governed by the stochastic differential equation
\begin{equation}\label{stoc}
\left\{
   \begin{array}{l}
       dX=\sigma^T P^{T}dW\\
    X(0)=\xi
    \end{array}
    \right.,
    \end{equation}
    where $W$ denotes the standard 
    $2n-$dimensional Brownian motion.\\
    Indeed, taking into account that $(X_iX_j\phi)_{sym}=\sigma (\xi)D^2\phi\sigma^T (\xi)$, it is immediate to verify that 
    $$\mbox{tr}(\sigma^{T}P^{T}P\sigma
D^{2}\phi)=\mbox{tr}(P^{T}P\sigma D^{2}\phi\sigma^{T})=\mbox{tr}(AD^2_{H^n}\phi).
$$
The operator $\tilde{\cal  M}^+_{\lambda,\Lambda}$ being  the supremum of $L_A$, is involved in controlled diffusion processes of the same type as (\ref{stoc}).

The first aim of this paper is to construct some functions related to the Pucci--Heisenberg operators $ \tilde{\cal M}^\pm_{\lambda,\Lambda}$  which will play the role of the ''fundamental solutions'' in this setting and which reduce  to the standard fundamental solution of the Heisenberg Laplacian $-\Delta_{H^n}$, found by Folland (see \cite{folland}), in the linear case $\lambda=\Lambda$.  These functions were constructed for the uniformly elliptic operators ${\cal  P}^{\mp}_{\lambda,\Lambda}$ by Pucci in \cite{puc} and called ''extremal barriers''.

The construction of these ''fundamental solutions''  for the Pucci--Heisenberg  operators $ \tilde{\cal M}^\pm_{\lambda,\Lambda}$, shall allow us  to answer to two different questions about the qualitative properties of solutions of (\ref{0}).\\
The first one concerns the existence and the uniqueness of the viscosity solution, continuous up to the 
boundary to Dirichlet problem associated with $F$ in a bounded, regular  domain $\Omega$, see Section \ref{dir}.\\
The second one is to state qualitative properties of Hadamard, Liouville and Harnack type for viscosity solutions associated with $F$, generalizing the analogues obtained for the uniformly elliptic fully nonlinear case in \cite{acfl, cdc2}; see Section  \ref{qua}. \\

Our existence results are based on  the Comparison Principle (see\cite{bardiman,wang}) and on the Perron--Ishii method (see \cite{user}) as in the paper by Bardi\&Mannucci, see \cite{bardiman}.\\
Their existence results,  in the case of fully non linear operators constructed with Heisenberg vector fields, are  only given for constant boundary data and for non linear operators obtained  by adding to the Pucci-Heisenberg operator $ \tilde{\cal M}^+_{\lambda,\Lambda}$  an operator depending on first order term which does not degenerate at the characteristic points of the
boundary (which is assumed to be either the Heisenberg ball, defined in Section \ref{pre}, or the euclidean one).\\
Precisely, they impose the following structure condition on the first order term:
\[
\begin{array}{l}
H(x,0)\leq  0\,,\\
H(x,p)\geq H_{hom}(x,p)-M\,,
\end{array}
\]
where $H_{hom}$ is a $1-$homogeneous function with respect to the gradient variable, which does not degenerate at the characteristic points of the boundary.\\
 These conditions allow the authors to use as lower barrier the boundary datum and as upper barrier the fundamental solution  for the Heisenberg Laplacian, see \cite{folland}. \\
On the other hand, their existence results apply to more general fully
nonlinear degenerate elliptic equations, not necessarily related to the
Heisenberg vector fields.\\
As far as the issues of existence and uniqueness in the fully non linear degenerate case are concerned, we mention the papers \cite{b, w2} and the references therein, where  Dirichlet problems for the infinity laplacian on Carnot groups were studied, and the work in progress by Birindelli\&Stroffolini,\cite{birstr} where some holder continuity properties of solutions to Dirichlet problems similar to (\ref{0}) are studied.

Our first existence result is relative to the following Dirichlet problem:
\begin{equation}\label{indiri}
\left\{
\begin{array}{ll}
F(\xi, D^2_{H^n}u)=0\qquad &\mbox{in }\Omega\\
u=\psi&\mbox{on }\partial\Omega\,,
\end{array}
\right .
\end{equation}
where  $\Omega$ is a bounded domain verifying the exterior Heisenberg--ball condition at the boundary (see Definition \ref{sferaext}), and  $\psi\in C(\partial\Omega)$. 

The second result involves operators depending on first order terms which have to degenerate at the characteristic points of the boundary, see (\ref{hgen}).\\
A preliminary result concerns the Dirichlet problem for the equation (\ref{0}) in the model annular domain  (see also \cite{crasalenic}).\\
Then, by using the barriers of the model problem, we construct local barriers for  the same Dirichlet problem  in a bounded regular domain satisfying the exterior Heisenberg--ball condition at the characteristic points of the boundary. Let us point out that this  geometrical condition on the domain is stronger than the classical ones (intrinsic cone condition, capacity,\dots) for the linear or semilinear  problems associated with the Heisenberg laplacian $\Delta_{H^n}$, see e.g. \cite{gaveau,gallardo,jerison,bcdc2}.\\

Such barriers allow us to prove the existence and uniqueness of the viscosity solution, continuous up to the boundary, for:
\begin{equation}\label{indiri2}
\left\{
\begin{array}{ll}
F(\xi, D^2_{H^n}u)+H(\xi,\nabla_{H^n}u)=0\qquad &\mbox{in }\Omega\\
u=\psi&\mbox{on }\partial\Omega\,.
\end{array}
\right .
\end{equation}
In Section \ref{qua} we first generalize some  Hadamard type theorems to this fully nonlinear degenerate setting, by proving that  the minima on annular sets of viscosity  supersolutions to the Pucci--Heisenberg operators, satisfy some ''concavity'' relations with respect to the associated ''fundamental solutions'' (see Theorem \ref{hadamdeg}).\\
Further,  we deduce  some non-linear Liouville results for supersolutions of $ \tilde{{\cal  M}}^{-}_{\lambda,\Lambda}$ (or subsolutions of $ \tilde{{\cal  M}}^{+}_{\lambda,\Lambda}$) in the whole space $\R^{2n+1}$ when the dimension $2n+1\leq \frac{\Lambda}{\lambda}$ , depending on the fact that in this case  the ''fundamental solutions'' blow up at infinity. \\
Finally, we state a weak Harnack inequality for radial supersolutions  with respect to $\rho$ (defined in  (\ref{normaodeltah})) of the inequality  $F(\xi,D^2_{H^n}u)\geq 0$.

\noindent{\bf Acknowledgment:} This paper was mainly done while the first author was  visiting the University of Rennes 1  with the support of the IRMAR. She wish to thank the people of the IRMAR for the kind hospitality.

\section{Preliminary facts}
\label{pre}
In this section we briefly recall some basic facts regarding  the Heisenberg vector fields $X_i$ defined in (\ref{heisf}) as well as the definition and some properties of viscosity solutions to fully nonlinear operators which will be useful in the sequel. For more details, see e.g. \cite{folland,folland2,hormander,user,caff,ipl}\\
First of all, let us point out that $X_i$ satisfy:
\[
[X_i,X_{i+n}]=-4\frac{\partial}{\partial\xi_{2n+1}},\qquad [X_i,X_j]=0\quad\forall
j\ne i+n; i\in\{1,\dots,n\}
\]
and therefore generate the whole Lie algebra  of left-invariant vector fields on the Heisenberg group 
 $H^n=(\R^{2n+1},\circ ) $ where $\circ$ is defined by 
\begin{equation}\label{circ}
\eta\circ\xi= 
(\xi_1+\eta_1,\dots,\xi_{2n}+\eta_{2n},\xi_{2n+1}+\eta_{2n+1}+
2 \sum_{i=1}^n (\xi_{i}\eta_{i+n}-\xi_{i+n}\eta_{i})).
\end{equation}
Moreover, they are homogeneous of degree $-1$
with respect to the {\it anistropic dilations}
\begin{equation}\label{dil}
\delta_\lambda (\xi)=
(\lambda\xi_{1},\dots,\lambda\xi_{2n},\lambda^{2}\xi_{2n+1})\;\;
(\lambda >0\;)\,.
\end{equation}
Then it is useful to consider the following
{\it homogeneous norm} which is $1-$homogeneous with respect to (\ref{dil})
\begin{equation}
\label{normaodeltah}
\rho(\xi)=(( \sum_{i=1}^{2n} \xi_i^2)^2  +
\xi_{2n+1}^{2})^{\frac{1}{4}}\;\;,
\end{equation} 
and the associated Heisenberg distance
\begin{equation}\label{dh}
d_H(\xi,\eta)=\rho(\eta^{-1}\circ\xi)\,.
\end{equation}
We denote by $B^H_R(\xi)$ the Koranyi-ball with centre at $\xi$ and radius $R$ associated with the distance (\ref{dh}), and  we will refer to it as {\sl  Heisenberg-ball.} The volume of $B^H_R(\xi)$ does not depend on the centre $\xi$ and scales as $R^Q$, where $Q=2n+2$ is the homogeneous dimension of the Heisenberg group.\\
Moreover, if we consider the $(2n)\times (2n+1)$ matrix whose rows are the coordinates of the vector fields $X_i$, that is 
\begin{equation}\label{sigmadeg}\sigma =\left (\begin{array}{ccc}
I_n&0&2y\nl
0&I_n&-2x
\end{array}\right ),
\end{equation}
where  $I_n$ is the identity $n\times
n$ matrix, $x=(\xi_1,\dots,\xi_n)^T$ and $y=(\xi_{n+1},\dots ,\xi_{2n})^T$ ($ ^T$ stands for transposition), then the Heisenberg gradient $\nabla_{H^n}$ of a function $\phi\,:\,\R^{2n+1}\rightarrow  \R$ is expressed by
\[\nabla_{H^n}\phi=(X_1\phi,\dots ,X_{2n}\phi)=\sigma (\xi)\nabla\phi\]
(where $\nabla$ stands for the usual gradient), and  the Heisenberg hessian of $\phi$ 
\[D^2_{H^n}\phi=(X_iX_j\phi)_{sym}=\sigma (\xi)D^2\phi\sigma^T (\xi),\]
where $sym$ denotes the symmetrized matrix.
\bigskip
 
As mentioned in the Introduction, we are interested in the fully non linear equation (\ref{0}), where the second order term is of the form $F(\xi,\sigma (\xi)M\sigma^T (\xi)) $ and satisfies,  for some
$0<\lambda\leq
\Lambda$, the following conditions:
\begin{equation}\label{sigmaue}
\lambda \, \mbox{tr}(\sigma (\xi)P\sigma^T(\xi))\leq
F(\xi,\sigma (\xi)(M)\sigma^T(\xi))-F(\xi,\sigma (\xi)(M+P)\sigma^T(\xi))\leq
\Lambda \, 
\mbox{tr}(\sigma (\xi)P\sigma^T(\xi))
\end{equation}
for all $M,P\in {\cal S}_{2n+1}$ with $P\geq
0$ (i.e. nonnegative definite) and 
\begin{equation}\label{sigmaf00} F(\xi ,0)\equiv 0\qquad\mbox{ for all }\xi\in\R^{2n+1}\, .
\end{equation}

The most important examples of second order operators verifying (\ref{sigmaue}) and (\ref{sigmaf00}) are the Pucci-Heisenberg operators (\ref{puc}). 
Moreover, being the vector fields $X_i$ invariant with respect to the left-translation (\ref{circ}), we get that 
$ \tilde{\cal  
M}^{\mp}_{\lambda,\Lambda} $ are invariant with respect the left-translations too.  

Furthermore, let us observe that the operators $  \tilde{\cal  
M}^-_{\lambda,\Lambda}(M)$ and $ \tilde{\cal  
M}^+_{\lambda,\Lambda}(M)$ can be rewritten respectively as 
$$ \tilde{\cal  M}^-_{\lambda,\Lambda}(M)= \inf_{\lambda I\leq B\leq\Lambda I }-\mbox{tr}(B \sigma(\xi)M\sigma (\xi)^T )\,,$$
$$ \tilde{\cal  M}^+_{\lambda,\Lambda}(M)= \sup_{\lambda I\leq B\leq\Lambda I }-\mbox{tr}(B \sigma(\xi)M\sigma (\xi)^T )\,.$$
Thus, setting $B=P^TP$, and using that $tr(P^TP\sigma(\xi)M\sigma(\xi)^{T})=tr((P\sigma(\xi))^TP\sigma(\xi)M)$, we get
$$ \tilde{\cal  M}^-_{\lambda,\Lambda}(M)=\inf_{P}-\mbox{tr}((P \sigma(\xi))^TP\sigma (\xi) M)$$
 and an analogous expression holds  for the operator $\tilde{\cal  M}^+_{\lambda,\Lambda}$.\\
Moreover, 
\begin{equation}\label{sigmaeqmpm}
 \tilde{\cal M}^+_{\lambda,\Lambda}(M)=- \tilde{\cal M}^-_{\lambda,\Lambda}(-M)\qquad\forall \  
M\in{\cal S}_{2n+1}.
\end{equation}

 Let us recall the definition of viscosity solution to 
\begin{equation}\label{equa}F(\xi, D^2_{H^n}u)+H(\xi,\nabla_{H^n}u)=0\qquad \mbox{in }\Omega \,.
\end{equation}

\begin{defi}\label{defsol} {\em
A lower semicontinuous function $v:\Omega\subseteq \R^{2n+1}\to \R$ is a {\sl viscosity
supersolution } of (\ref{equa}) 
if, for all $\zeta \in C^2(\Omega )$ and  $\xi_0\in\Omega $ such that $v-\zeta $ has a
local minimum at $\xi_0$, it results
\[F(\xi_0,D^2_{H^n}\zeta (\xi_0) )+H(\xi_0,\nabla_{H^n}\zeta (\xi_0) )\geq  0\,
\]
whereas,  an upper semicontinuous function $u:\Omega\subseteq \R^{2n+1}\to \R$ is a viscosity subsolution of (\ref{equa}) if for all $\zeta \in C^2(\Omega )$ and  $\eta_0\in\Omega $ such that $u-\zeta $ has a
local maximum at $\eta_0$, it results
\[F(\eta_0,D^2_{H^n}\zeta (\eta_0) )+H(\eta_0,\nabla_{H^n}\zeta (\eta_0) )\leq  0\,.
\]
Moreover we say that $w\in C(\Omega)$ is a viscosity solution of (\ref{equa}) if it simultaneously is a viscosity sub and supersolution.}

\end{defi}

Since we are concerned with Dirichlet problems associated with (\ref{0}) and we shall investigate the continuity up to the boundary of viscosity solutions, we need to analyze the behaviour of these solutions at the characteristic points. Thus, let us recall their definition:
\begin{defi}\label{charpoint}
Let $\Omega$ be a bounded $C^2$ domain of $\R^{2n+1}$ such that 
\begin{equation}\label{regolarephi2}
\begin{array}{c}
\mbox{there exists }\Phi\in C^2(\overline{\Omega}) \mbox{ with }\nabla\Phi(\xi)\neq 0\quad\forall \xi\in\partial\Omega \mbox{ such that :}\\
\Omega=\{\xi\in \R^{2n+1}\,:\,\Phi(\xi)>0\}
\\
\partial \Omega=\{\xi\in \R^{2n+1}\,:\,\Phi(\xi)=0\}\,.
\end{array}
\end{equation}
Then $\xi_0\in\partial\Omega$ is called a {\it characteristic} point of $\partial\Omega$ if:
\begin{equation}\label{defcar}
\nabla_{H^n}\Phi(\xi_0)= 0\,.
\end{equation}
\end{defi}
Let us observe that, if $\xi_0$ is a characteristic point, then necessarily $\frac{\partial \Phi}{\partial\xi_{2n+1}}(\xi_0)\ne 0$ and the notion of being characteristic  does not depend on the choice of $\Phi$.\\

Finally, let us state  the following version of Comparison Principle recently proved by Bardi \& Mannucci (see \cite{bardiman}), which shall allow us to use the Perron-Ishii method in order to get the existence and uniqueness of continuous (up to the boundary) viscosity solutions to Dirichlet problems for  (\ref{0}).
\begin{thm}\label{cp}
Let $F$ satisfy (\ref{sigmaue}) and (\ref{sigmaf00}) and $H$ satisfy (\ref{inH}); let $u\in BUSC (\Omega)$ be a viscosity subsolution of $F+H$ and $v\in BLSC (\Omega)$ be a viscosity supersolution of $F+H$ (where $BUSC$ and $BLSC$ respectively stand for bounded upper semicontinuous and bounded lower semicontinuous function). If $u\leq v$ on $\partial\Omega$, then $u\leq v$ in $\overline\Omega$.
\end{thm}
\section{''Fundamental solutions''  for Pucci-Heisenberg operators} \label{fon}
In this section  we will construct some radial functions for the operators $\tilde{\cal  M}^+_{\lambda,\Lambda}$ and $\tilde{\cal  M}^-_{\lambda,\Lambda}$ wich will play the role of the ''fundamental solution'' in this nonlinear framework.\\
 To this purpose, we will make use of the expression of $ \tilde{\cal 
 M}^+_{\lambda,\Lambda}$ on
functions $f$ depending on the homogeneous norm defined in (\ref{normaodeltah}). 
\subsection{Pucci--Heisenberg operators on radial functions}
Let us observe that, for $\rho (\xi)>0$,
\begin{equation}\label{radialhessprov}
   \displaystyle
D^2_{H^n}f(\rho)=f'(\rho)D^2_{H^n}\rho+f^{\prime\prime}(\rho)\nabla_{H^n}\rho\otimes\nabla_{H^n}\rho\,;
    \end{equation}
moreover, the following technical lemma, provides the expression of the 
intrinsic Heisenberg hessian matrix $D^{2}_{H^n}\rho$.

\begin{lem}\label{d2rho}
    Let $\rho(\xi)$ be defined in (\ref{normaodeltah}), then
\[D^2_{H^n}\rho=-\frac{3}{\rho}\nabla_{H^n}\rho\otimes\nabla_{H^n}\rho+
\frac{1}{\rho}|\nabla_{H^n}\rho|^2I_{2n}+\frac{2}{\rho^3}\left(\begin{array}{cc}
B&C\nl
-C&B
\end{array}
\right)
\]
where the matrices $B=(b_{i,j})$ and $C=(c_{i,j})$ are defined as follows:
\[b_{i,j}=\xi_i\xi_j+\xi_{n+i}\xi_{n+j}\qquad c_{i,j}=\xi_i\xi_{n+j}-\xi_j\xi_{n+i},
\]
for $ i,j\in\{1,\dots ,n\}$ (observe that $B$ is symmetric whereas 
$C^{T}=-C$),  and 
\[ |\nabla_{H^n}\rho|^2=\frac{\sum_{i=1}^{2n} \xi_i^2}{\rho^2}\qquad\mbox{ for
}\rho>0.
\]
\end{lem}

\noindent {\it Proof.}
A standard calculation shows that for $i=1,\dots,n$:
\[\displaystyle\begin{array}{l}
X_{i}\rho=\xi_{i}\frac{|\nabla_{H^n}\rho|^2}{\rho}+\frac{\xi_{i+n}\xi_{2n+1}}{\rho^{3}}\\
X_{n+i}\rho=\xi_{i+n}\frac{|\nabla_{H^n}\rho|^2}{\rho}-\frac{\xi_{i}\xi_{2n+1}}{\rho^{3}}\,.
\end{array}
\]
Moreover, for $i,j=1,\dots ,n\quad i\ne j$, we have:
\begin{equation}\label{hes1}\displaystyle
\begin{array}{l}
X_{j}X_{i}\rho=\frac{2}{\rho^{3}}b_{ij}-\frac{3}{\rho^{7}}\left[\xi_{i}\rho^{2}|\nabla_{H^n}\rho|^2+\xi_{i+n}\xi_{2n+1}\right]
\left[\xi_{j}\rho^{2}|\nabla_{H^n}\rho|^2+\xi_{j+n}\xi_{2n+1}\right]\\
X_{n+j}X_{n+i}\rho=\frac{2}{\rho^{3}}b_{ij}-\frac{3}{\rho^{7}}\left[\xi_{n+i}\rho^{2}|\nabla_{H^n}\rho|^2-\xi_{i}\xi_{2n+1}\right]
\left[\xi_{n+j}\rho^{2}|\nabla_{H^n}\rho|^2-\xi_{j}\xi_{2n+1}\right]\\
X_{n+j}X_{i}\rho=\frac{2}{\rho^{3}}c_{ij}-\frac{3}{\rho^{7}}\left[\xi_{i}\rho^{2}|\nabla_{H^n}\rho|^2+\xi_{i+n}\xi_{2n+1}\right]
\left[\xi_{n+j}\rho^{2}|\nabla_{H^n}\rho|^2-\xi_{j}\xi_{2n+1}\right]
\end{array}
\end{equation}
and
\begin{equation}\label{hes2}\displaystyle\begin{array}{l}
X_{i}^{2}\rho=\frac{|\nabla_{H^n}\rho|^2}{\rho}+\frac{2}{\rho^{3}}b_{ii}-
\frac{3}{\rho^{7}}\left[\xi_{i}\rho^{2}|\nabla_{H^n}\rho|^2+\xi_{i+n}\xi_{2n+1}\right]^{2}\\
X_{i+n}^{2}\rho=\frac{|\nabla_{H^n}\rho|^2}{\rho}+\frac{2}{\rho^{3}}b_{ii}-
\frac{3}{\rho^{7}}\left[\xi_{i+n}\rho^{2}|\nabla_{H^n}\rho|^2-\xi_{i}\xi_{2n+1}\right]^{2}\\
X_{n+i}X_{i}\rho=\frac{\xi_{{2n+1}}}{\rho^{3}}-\frac{3}{\rho^{7}}\left[\xi_{i}\rho^{2}|\nabla_{H^n}\rho|^2+\xi_{i+n}\xi_{2n+1}\right]
\left[\xi_{n+i}\rho^{2}|\nabla_{H^n}\rho|^2-\xi_{i}\xi_{2n+1}\right]\\
X_{i}X_{n+i}\rho=-\frac{\xi_{{2n+1}}}{\rho^{3}}-\frac{3}{\rho^{7}}\left[\xi_{i}\rho^{2}|\nabla_{H^n}\rho|^2+\xi_{i+n}\xi_{2n+1}\right]
\left[\xi_{n+i}\rho^{2}|\nabla_{H^n}\rho|^2-\xi_{i}\xi_{2n+1}\right]\,.
\end{array}
\end{equation}
Now, taking into account that
\begin{equation}\label{tensor}
    \displaystyle\begin{array}{l}
X_{i}\rho\,X_{j}\rho=\frac{1}{\rho^{6}}\left[\xi_{i}\rho^{2}|\nabla_{H^n}\rho|^2+\xi_{i+n}\xi_{2n+1}\right]
\left[\xi_{j}\rho^{2}|\nabla_{H^n}\rho|^2+\xi_{j+n}\xi_{2n+1}\right]\\
X_{n+i}\rho\,X_{n+j}\rho=\frac{1}{\rho^{6}}\left[\xi_{n+i}\rho^{2}|\nabla_{H^n}\rho|^2-\xi_{i}\xi_{2n+1}\right]
\left[\xi_{n+j}\rho^{2}|\nabla_{H^n}\rho|^2-\xi_{j}\xi_{2n+1}\right]\\
X_{i}\rho\,X_{n+j}\rho=\frac{1}{\rho^{6}}\left[\xi_{i}\rho^{2}|\nabla_{H^n}\rho|^2+\xi_{i+n}\xi_{2n+1}\right]
\left[\xi_{n+j}\rho^{2}|\nabla_{H^n}\rho|^2-\xi_{j}\xi_{2n+1}\right]\,
\end{array}
\end{equation}
by substituting (\ref{tensor}) in (\ref{hes1}) and (\ref{hes2}), we get the claim\,.$\hfill\Box$

 \par

\medskip
By using the expression of $D^2_{H^n}\rho$ provided in Lemma 
\ref{d2rho}, (\ref{radialhessprov}) becomes:
\begin{equation}\label{radialhess}
   \displaystyle\begin{array}{ll}
    D^2_{H^n}f(\rho)=&\frac{f'(\rho)}{\rho}\left [ 
    -3\nabla_{H^n}\rho\otimes\nabla_{H^n}\rho + 
    |\nabla_{H^n}\rho|^2I_{2n}+\frac{2}{\rho^2}\left(\begin{array}{cc}
B&C\nl
-C&B
\end{array}
\right) \right]\\
&+f^{\prime\prime}(\rho)\nabla_{H^n}\rho\otimes\nabla_{H^n}\rho\,.
 \end{array}
\end{equation}
    
\begin{lem}\label{lemmaeigen}
The eigenvalues of $D^2_{H^n}f(\rho)$ are 
\begin{equation}\label{eigenvhn}
\begin{array}l
|\nabla_{H^n}\rho|^2f^{\prime\prime}(\rho)\mbox{ which is simple  }\\
3|\nabla_{H^n}\rho|^2\frac{f'(\rho)}{\rho}\mbox{ which is simple}\\
|\nabla_{H^n}\rho|^2\frac{f'(\rho)}{\rho}\mbox{  which has multiplicity
$2n-2$,}
\end{array}
\end{equation}
if $\nabla_{H^n}\rho\ne 0$; otherwise all the eigenvalues of $D^2_{H^n}f(\rho)$ 
vanish identically. 
\end{lem}

\noindent {\it Proof.} 
Assume that $\nabla_{H^n}\rho\ne 0$, the other case being obvious.
Then, the vector $w=\frac{\nabla_{H^n}\rho}{|\nabla_{H^n}\rho|}$ is an
eigenvector associated with 
$|\nabla_{H^n}\rho|^2f^{\prime\prime}(\rho)$. Indeed, being
$p\otimes p\frac{p}{|p|}=p|p|$, and
\[
\left(\begin{array}{cc}
B&C\nl
-C&B
\end{array}
\right)\nabla_{H^n}\rho=|\nabla_{H^n}\rho|^2\rho^{2}\nabla_{H^n}\rho
\]
we get that
\[\begin{array}{l}
D^2_{H^n}f(\rho)w=f^{\prime\prime}(\rho)|\nabla_{H^n}\rho|^{2}w\,.
 \end{array}
 \]
Hence, $|\nabla_{H^n}\rho|^2f^{\prime\prime}(\rho)$ is an eigenvalue.

Now, let us consider the vector 
$v=(X_{n+1}\rho,X_{n+2}\rho,\dots,X_{2n}\rho,-X_{1}\rho,-X_{2}\rho,\dots ,-X_{n}\rho)$.
Since $v\perp \nabla_{H^n}\rho$, 
\[
D^2_{H^n}f(\rho)v=\frac{f^{\prime}(\rho)}{\rho}\left
[|\nabla_{H^n}\rho|^{2}v+\frac{2}{\rho^{2}}\left(\begin{array}{cc}
B&C\nl
-C&B
\end{array}
\right)v\right]\,,
 \]
and moreover,
\[\left(\begin{array}{cc}
B&C\nl
-C&B
\end{array}
\right)v=|\nabla_{H^n}\rho|^{2}\rho^{2}v\,.
\]
Hence, 
\[\displaystyle 
D^2_{H^n}f(\rho)v=3\frac{f^{\prime}(\rho)}{\rho}|\nabla_{H^n}\rho|^{2}v\,,
\]
which proves that
$3\frac{f^{\prime}(\rho)}{\rho}|\nabla_{H^n}\rho|^{2}$ is another
eigenvalue.\\
Now, in order to prove that
$|\nabla_{H^n}\rho|^2\frac{f'(\rho)}{\rho}$ is an eigenvalue with multiplicity
$2n-2$, we need to find $2n-2$ eigenvectors which are orthogonal to
$\{v,w\}$ and belong to the $\mbox{ Ker }\left(\begin{array}{cc}
B&C\nl
-C&B
\end{array}\right)
$. \\
At this purpose, let us choose, for
$k=1,\dots\,,\left[\frac{n}{2}\right ] $, the vectors
$\eta=(\eta_{1},\dots\,,\eta_{2n})$
defined as 
\[\displaystyle\begin{array}{ll}
\eta_{2k}&=X_{2k}\rho X_{n+2k-1}\rho-X_{2k-1}\rho X_{n+2k}\rho\\
\eta_{n+2k-1}&=-(X_{2k}^{2}\rho + X_{n+2k}^{2}\rho)\\
\eta_{n+2k}&=X_{2k-1}\rho X_{2k}\rho+X_{n+2k-1}\rho X_{n+2k}\rho\\
\eta_{j}&=0\qquad\mbox{ for }j\ne 2k\,,n+2k-1\,,n+2k\,.
\end{array}
\]
It is immediate to verify that $\eta\perp v$ and
$\eta\perp w$. Moreover,  being:
\[\displaystyle{\begin{array}{ll}
\eta_{2k}&=\left(\frac{|\nabla_{H^n}\rho|^4}{\rho^2}+\frac{\xi_{2n+1}^2}{\rho^4}\right )c_{2k,2k-1}\\
\eta_{n+2k-1}&=-\left(\frac{|\nabla_{H^n}\rho|^4}{\rho^2}+\frac{\xi_{2n+1}^2}{\rho^4}\right )b_{2k,2k}\\
\eta_{n+2k}&=\left(\frac{|\nabla_{H^n}\rho|^4}{\rho^2}+\frac{\xi_{2n+1}^2}{\rho^4}\right )b_{2k,2k-1}\\
\eta_{j}&=0\qquad\mbox{ for }j\ne 2k\,,n+2k-1\,,n+2k\,,
\end{array}}
\]
and taking into account that:
\[\displaystyle{\begin{array}{ll}
b_{i,2k}c_{2k,2k-1}-b_{2k,2k}c_{i,2k-1}+b_{2k,2k-1}c_{i,2k}=0\qquad &\mbox{ for }i=1,\dots,n\,,\\
b_{i,2k}b_{2k,2k-1}-c_{i,2k}c_{2k,2k-1}-b_{i,2k-1}b_{2k,2k}=0&\mbox{ for }i=1,\dots,n\,,
\end{array}}
\]
an easy calculation shows that
$$\eta\in\mbox{ Ker }\left(\begin{array}{cc}
B&C\nl
-C&B
\end{array}\right)\,.
$$
Hence, the $\left[\frac{n}{2}\right ] $ vectors $\eta$ are
eigenvectors associated with $|\nabla_{H^n}\rho|^2\frac{f'(\rho)}{\rho}$. Other $\left[\frac{n}{2}\right ] $ eigenvectors related to the same eigenvalue are simply obtained by setting:
\[\displaystyle\begin{array}{ll}
\hat{\eta}_{2k-1}&=\eta_{n+2k-1}\\
\hat{\eta}_{2k}&=\eta_{n+2k}\\
\hat{\eta}_{n+2k}&=-\eta_{2k}\\
\hat{\eta}_{j}&=0\qquad\mbox{ for }j\ne 2k-1\,,2k\,,n+2k\,.
\end{array}
\]
If $n$ is odd,  other two eigenvectors associated with the same eigenvalue 
are respectively,  the vector $\eta$ defined by:
\[\displaystyle
{\begin{array}{ll}
\eta_{n}&=-\frac{c_{n,n-1}}{b_{n,n}}=
\frac{\xi_{n-1}\xi_{2n}-\xi_{n}\xi_{2n-1}}{\xi_{n}^{2}+\xi_{2n}^{2}}\\
\eta_{2n-1}&=1\\
\eta_{2n}&=-\frac{b_{n-1,n}}{b_{n,n}}=
-\frac{\xi_{n-1}\xi_{n}+\xi_{2n}\xi_{2n-1}}{\xi_{n}^{2}+\xi_{2n}^{2}}\\
\eta_{j}&=0\qquad\mbox{ if }j\ne n\,,2n-1\,,2n\,,
\end{array}}
\]
and $\hat{\eta}$ defined by:
\[\displaystyle\begin{array}{ll}
\hat{\eta}_{n-1}&=\eta_{2n-1}\\
\hat{\eta}_{n}&=\eta_{2n}\\
\hat{\eta}_{2n}&=-\eta_{n}\\
\hat{\eta}_{j}&=0\qquad\mbox{ for }j\ne n-1\,,n\,,2n\, ,
\end{array}
\]
as it can be seen by using that:
\[\displaystyle{\begin{array}{ll}
-b_{i,n}c_{n,n-1}+b_{n,n}c_{i,n-1}-b_{n-1,n}c_{i,n}=0\qquad &\mbox{ for }i=1,\dots,n\,,\\
-b_{i,n}b_{n-1,n}+c_{i,n}c_{n,n-1}+b_{i,n-1}b_{n,n}=0&\mbox{ for }i=1,\dots,n\,.
\end{array}}
\]
The other $n-2$ eigenvectors can be chosen of the form
$(\eta_{1},\eta_{2}\,,\eta_{n},0\,,0\,\dots\,,0)$, with 
$$\displaystyle\begin{array}{l}
(\eta_{1}\,,\eta_{2}\,,\dots,\eta_{n})\perp (\xi_{1}\,,\xi_{2}\,,\dots
,\xi_{n})\\
(\eta_{1}\,,\eta_{2}\,,\dots,\eta_{n})\perp (\xi_{n+1}\,,\xi_{n+2}\,,\dots
,\xi_{2n})\,.
\end{array}
$$
Indeed, setting $\tilde\eta=(\eta_{1}\,,\eta_{2}\,,\dots,\eta_{n})$, 
we get
\[\displaystyle\begin{array}{l}
B\tilde\eta\cdot\tilde\eta=\left(\sum_{i=1}^{n}\xi_{i}\eta_{i}\right)^{2}+
\left(\sum_{i=1}^{n}\xi_{n+i}\eta_{i}\right)^{2}=0\\
C\tilde\eta\cdot\tilde\eta=\left(\sum_{i=1}^{n}\xi_{i}\eta_{i}\right)
\left(\sum_{j=1}^{n}\xi_{n+j}\eta_{j}\right)-\left(\sum_{j=1}^{n}\xi_{j}\eta_{j}\right)
\left(\sum_{i=1}^{n}\xi_{n+i}\eta_{i}\right)=0\,.
\end{array}
\]
Moreover,
\[\displaystyle \tilde\eta\perp (X_{1}\rho (\xi),\dots X_{n}\rho
(\xi))\qquad \mbox{ and }\tilde\eta\perp (X_{n+1}\rho (\xi),\dots X_{2n}\rho
(\xi))\,.
\]
Hence we easily conclude that the multiplicity of $|\nabla_{H^n}\rho|^2\frac{f'(\rho)}{\rho}$ 
is $2n-2$. This concludes the proof of the lemma\,.$\hfill\Box$

 \par
Lemma \ref{lemmaeigen}
 implies immediately the following

\begin{cor}\label{radialcor}
For every function $\Phi(\xi)=\varphi (\rho(\xi))$ with $\varphi $ concave and increasing, it results:
\begin{equation}\label{sigmam-ci}
\tilde{\cal M}^+_{\lambda ,\Lambda }(D^2 \Phi )=
|\nabla_{H^n}\rho|^2[-\lambda (2n+1)\frac{\varphi ^\prime (\rho)}{\rho} - \Lambda \varphi ^{\prime
\prime} (\rho)]
\end{equation}
 whereas if $\varphi $  is  convex and decreasing, then
\begin{equation}\label{sigmam-cd}
\tilde{\cal M}^+_{\lambda ,\Lambda }(D^2 \Phi )=
|\nabla_{H^n}\rho|^2[-\lambda \varphi ^{\prime \prime}(\rho) - \Lambda (2n+1)
\frac{\varphi ^\prime
(\rho)}{\rho}]\,.
\end{equation}
$\hfill\Box$
\end{cor}

\subsection{''Fundamental solutions''}

Using (\ref{sigmam-ci}), and (\ref{sigmam-cd}), we can
construct  radial functions with respect to (\ref{normaodeltah}) $\Phi
(\xi)=\varphi (\rho(\xi))$ which are classical solutions of the equations:
\begin{equation}\label{sigmaomopm}
\tilde{\cal M}^\pm_{\lambda ,\Lambda }(D^2 \Phi )=0\qquad \mbox{ in\ }\  \R^{2n+1}
\setminus \{0\} \, ,
\end{equation}
and are either concave and increasing or convex and
decreasing.\\
These functions  play the same role as the fundamental
solution in many qualitative properties that will be detailed in the next sections (see the book \cite{pw} for the classical harmonic functions). It is in this respect that
they will be considered as the \lq\lq fundamental solutions" of the equations (\ref{sigmaomopm}), as in \cite{puc}.
Let us point out, moreover, that in the particular case in which $\lambda
=\Lambda$, equations (\ref{sigmaomopm}) reduce to the Heisenberg Laplace equation; in this case, they coincide with  the classical fundamental solution for the Heisenberg Laplacian found by Folland in \cite{folland}.
These functions depend on the parameters $\alpha\geq 1$ and $\beta\geq 4$ defined by:
\begin{equation}\label{sigmanewN1}
\alpha \, =\, \frac{\lambda}{\Lambda} (Q-1) +1
\end{equation}
\begin{equation}\label{sigmanewN2}
\beta \, =\, \frac{\Lambda}{\lambda} (Q-1) +1
\end{equation}
($Q=2n+2$) which have to be considered as new dimensions, depending on the nonlinearities, and which coincide with the intrinsic linear dimension Q  if and only if $\lambda=\Lambda$.
\begin{lem}\label{funsolm+}
The radial functions
\[
\Phi_1 (\xi)\, =\, \varphi_1 (\rho(\xi))
\qquad\mbox{   and   }\qquad
\Phi_2 (\xi)\, =\, \varphi_2 (\rho(\xi))\, ,
\]
with $\varphi_1$ and $\varphi_2$ respectively given by 
\begin{equation}\label{sigmasolfon1}
\varphi_1 (\rho)=
\left\{ 
\begin{array}{ll}
C_1 \rho^{2-\alpha}+C_2 &\mbox{\ \ \ \  if \ \ } \alpha <2\nl
C_1\log \rho +C_2  &\mbox{\ \ \ \ if \ \ } \alpha =2\nl
-C_1 \rho^{2-\alpha}+C_2  &\mbox{\ \ \ \  if \ \ }
\alpha >2\, , \end{array}
\right.
\end{equation}
with constants $C_1\geq 0$ and $C_2\in \R$, 
and 
\begin{equation}\label{sigmasolfon2}
\varphi_2 (\rho)=C_1 \rho^{2-\beta }+C_2  
\end{equation}
are classical solutions (in particular, viscosity solutions) of
the equation:
\begin{equation}\label{sigmaomo+}
\tilde{\cal M}^+_{\lambda ,\Lambda }(D^2 \Phi )=0\qquad \mbox{ in\ }\  \R^{2n+1}
\setminus \{0\} \, .
\end{equation}
Moreover, $\varphi_1$ is concave and increasing whereas $\varphi_2$  is convex and
decreasing. 
\end{lem}
\noindent {\bf Proof}\quad                                                           
The statement immediately follows from  (\ref{sigmam-ci}) and (\ref{sigmam-cd}). Indeed,
the concave and increasing functions have to be looked for among the solutions of the ordinary
differential equation
\[
\lambda (2n+1)\frac{\varphi ^\prime (\rho)}{\rho} + \Lambda \varphi ^{\prime \prime}
(\rho) =0\qquad \mbox{  in\ }\  (0,+\infty )\, ,
\]
as well as the convex and decreasing solutions $\varphi$ must satisfy
\[
\lambda \varphi ^{\prime \prime}(\rho) + \Lambda (2n+1)
\frac{\varphi ^\prime
(\rho)}{\rho}
   =0\qquad \mbox{  in\ }\  (0,+\infty )\, .
\]
Then by an easy calculation, the claim follows.  $\hfill\Box$

\begin{rem}\label{remsolfonm-}
Remembering the relationship (\ref{sigmaeqmpm}) between $\tilde{\cal M}^+_{\lambda ,\Lambda}$
and $\tilde{\cal M}^-_{\lambda ,\Lambda}$, we have also found that the functions
\[
\Psi_1 (\xi)\, =\, -\Phi_2 (\xi)
\]
and
\[
\Psi_2 (\xi)\, =\, -\Phi_1 (\xi)
\]
are the \lq\lq fundamental solutions" of the equation
\begin{equation}\label{sigmaomo-}
\tilde{\cal M}^-_{\lambda ,\Lambda}(D^2 \Psi )=0\qquad \mbox{  in \  }\ \R^{2n+1}
\setminus \{0\} \, ,
\end{equation}
with $\Psi_1 (\xi)\equiv \psi_1 (\rho(\xi))=-\varphi_2 (\rho(\xi))$ such that $\psi_1 $ is a
concave and increasing function with respect to  $\rho(\xi)\in (0,\infty )$, and 
 $\Psi_2 (\xi)\equiv \psi_2 (\rho(\xi))=-\varphi_1 (\rho(\xi))$ such that $\psi_2 $ is a
convex and decreasing function with respect to $\rho$.
\end{rem}

If $\lambda=\Lambda$, we get  $\alpha =\beta =Q$ and the function $\Phi_1 \equiv \Phi_2$ coincides with the classical fundamental solution for the Heisenberg laplacian, see \cite{folland}.

Since the operators $\tilde{\cal M}^\pm_{\lambda ,\Lambda}$ are invariant with respect to the action (\ref{circ}) of the group, the  \lq\lq fundamental solutions"  which are singular at the point $\eta_0$ are  simply obtained  by substituting $\rho(\xi)$ with $\rho(\eta_0^{-1}\circ \xi)$.

The functions just constructed can be used as barriers in order to find domains which are regular for the Dirichlet problems. This will be the object of the following section.

\section{Existence of viscosity solutions to Dirichlet problems continuous up to the boundary}\label{dir}
\subsection{First case: Purely second order type operators}

We prove the existence of a continuous viscosity solution to the Dirichlet problem
\begin{equation}\label{diri}
\left\{
\begin{array}{ll}
F(\xi, D^2_{H^n}u)=0\qquad &\mbox{in }\Omega\\
u=\psi&\mbox{on }\partial\Omega\,,
\end{array}
\right .
\end{equation}
by using the Perron--Ishii method. 
Since the Comparison Principle holds true for (\ref{diri}) (see Theorem \ref{cp}),  it is enough to find, as it is well known, (see e.g. \cite{user}), a lower and an upper barrier (that is a subsolution and a supersolution of the equation in (\ref{diri}) which satisfy the boundary condition).

 The existence of such barriers is guaranteed if $\Omega$ satisfies the  {\it exterior Heisenberg--ball condition}. As  remarked in the introduction, this geometrical condition is stronger than the classical ones (intrinsic cone condition, capacity,\dots) in the linear variational case, see e.g. \cite{gaveau,gallardo}.
\begin{defi}\label{sferaext} We say that $\Omega$ satisfies the  exterior Heisenberg ball condition at $\xi_0\in\partial\Omega$ if \begin{equation}\label{sferapt}
\begin{array}{l}
\mbox{ there exist $\eta_0\in\Omega^C$, $r_0>0$ such that  $\overline{B^H_{r_0}(\eta_0)}\cap \overline{\Omega}=\xi_0$}\,.
\end{array}
\end{equation}
\end{defi}
\begin{lem}\label{barrierdiri}
Let $\Omega\subset\R^{2n+1}$ be a bounded domain satisfying the  {\it exterior  Heisenberg--ball condition} at all points of $\partial\Omega$, let $F$ verify (\ref{sigmaue}) and (\ref{sigmaf00})
and let $\psi$ be a continuous  function on $\partial\Omega$. Then, there exist a lower barrier $\underline u$ and an upper barrier $\overline{u}$ for the problem (\ref{diri}) which satisfy: $\underline{u}(\xi)=\overline{u}(\xi)=\psi(\xi)$ on $\partial\Omega$.
 \end{lem}

\noindent {\bf Proof}\quad
Let $\xi_0$ be a point of $\partial\Omega$.  
In order to construct  a local upper barrier for $F$ at $\xi_0$  it is enough to take a supersolution $v_{{\xi}_0}$ for  $\tilde{\cal M}^-_{\lambda ,\Lambda}$ such that $v_{{\xi}_0}(\xi_0)=0$ and $v_{{\xi}_0}>0$ in $\overline{\Omega}\setminus \{\xi_0\}$. Hence, let $v_{{\xi}_0}$ be the  solution of 
$$\begin{array}{ll}
\tilde{\cal M}^-_{\lambda ,\Lambda}(D^2v_{{\xi}_0})=0\qquad &\mbox{in }\R^{2n+1}\setminus\{\eta_0\}\\
v_{{\xi}_0}(\xi_0)=0\,,& 
\end{array}
$$
which is concave and increasing with respect to $\rho(\eta_0^{-1}\circ \xi)=d_{H^n}(\xi,\eta_0)$, that is 
$$v_{{\xi}_0}(\xi)=-\varphi_2 (\rho(\eta_0^{-1}\circ \xi))=r_0^{2-\beta}-d_{H^n}^{2-\beta}(\xi,\eta_0)\,,$$
where $\eta_0$ and $r_0$ are respectively the centre and the radius of the exterior Heisenberg--ball at $\xi_0$, defined in Definition \ref{sferaext}. 

Since $\psi$ is continuous at $ \xi_0 $, for every $\eps>0$ we can find $\delta_\eps>0$ such that  
$$|\psi(\xi)-\psi(\xi_0)|<\frac{\eps}{2}\qquad \mbox{ if }\xi\in\partial\Omega \mbox{ and }d_{H^n}(\xi,\xi_0)<\delta_\eps\,.$$
Let then $M=\sup_{\partial\Omega}|\psi|$ and $k>0$ be such that $k v_{{\xi}_0}(\xi)\geq 2M$ if $\xi\in\partial\Omega$ and $d_{H^n}(\xi,\xi_0)\geq\delta_\eps\,,$ and consider 

$$\overline{f}_{\eps,\xi_0}(\xi)=\psi(\xi_0)+\eps+k v_{{\xi}_0}(\xi)\qquad\mbox{ for } \xi\in\overline\Omega\,.$$
Then, 
$$0=\tilde{\cal M}^-_{\lambda ,\Lambda}(D^2 \overline{f}_{\eps,\xi_0}(\xi))\leq F(\xi,D^2_{H^n}\overline{f}_{\eps,\xi_0}(\xi))$$
and, moreover, $\overline{f}_{\eps,\xi_0}(\xi)\geq \psi (\xi)$ on $\partial\Omega$. 

Then, let us set
\begin{equation}\label{superbar}
\overline{u}(\xi):=\inf\{\overline{f}_{\eps,\xi_0}(\xi)\,\:\,\eps>0, \,\xi_0\in\partial\Omega\,\}\mbox{   for }\xi\in\overline\Omega\,.
\end{equation}
It is immediate to verify that $\overline{u}$ is lower semicontinuous and that it is  a viscosity supersolution of the equation (\ref{diri}). 

Moreover, for all $\xi\in\partial\Omega$, $\psi(\xi)\leq \overline{u}(\xi)\leq \overline{f}_{\eps,\xi}(\xi)=\psi(\xi)+\eps$, for all $\eps>0$. Hence, $\overline{u}=\psi$ on $\partial\Omega$ and it is an upper barrier for (\ref{diri}).

At the same way, using (\ref{sigmaeqmpm}), we can construct a lower barrier. Indeed, we easily have that
$$\underline{f}_{\eps,\xi_0}(\xi)=\psi(\xi_0)-\eps -k v_{{\xi}_0}(\xi)\qquad\mbox{ for } \xi\in\overline\Omega\,,$$
satisfies:
$$0=\tilde{\cal M}^+_{\lambda ,\Lambda}(D^2 \underline{f}_{\eps,\xi_0}(\xi))\geq F(\xi,D^2_{H^n}\underline{f}_{\eps,\xi_0}(\xi))$$
and 
$\underline{f}_{\eps,\xi_0}(\xi)\leq \psi (\xi)$ on $\partial\Omega$.

Thus,  setting
\begin{equation}\label{subbar}
\underline{u}(\xi):=\sup\{\underline{f}_{\eps,\xi_0}(\xi)\,\:\,\eps>0, \,\xi_0\in\partial\Omega\,\}\mbox{   for }\xi\in\overline\Omega\,,
\end{equation}
we get that $\underline{u}$ is upper semicontinuous and that it is  a viscosity subsolution of the equation (\ref{diri}).

Moreover, it achieves the boundary datum $\psi$  on $\partial\Omega$, as it is immediate to verify by the same argument as for $\overline{u}$. This concludes the proof.$\hfill\Box$

The previous lemma allows to apply the Perron-Ishii method in order to find a continuous viscosity solution of the problem (\ref{diri}). Indeed we have:
\begin{thm}\label{existence1}
 Let $\Omega\subset\R^{2n+1}$ be a bounded domain satisfying the  {\it exterior  Heisenberg--ball condition} at all points of $\partial\Omega$, let $F$ verify (\ref{sigmaue}) and (\ref{sigmaf00})
and let $\psi$ be a continuous  function on $\partial\Omega$. Then the function $W:\overline\Omega\rightarrow \R$ defined by:
 $$W(\xi)=\sup\{w(\xi)\,:\, \mbox{ $w$ is a subsolution of $F$ and }\underline{u}(\xi)\leq w(\xi)\leq \overline{u}(\xi) \mbox{ in }  \overline\Omega\}
 $$
 where $\underline{u}$ and $\overline{u}$ are defined in Lemma \ref{barrierdiri},  is the unique  continuous viscosity solution of (\ref{diri}) in $\overline\Omega$.
 \end{thm}
{\bf Proof.} In view of the Perron--Ishii method, using Lemma \ref{barrierdiri},  we get that the upper semicontinuous envelope of $W$ is a viscosity subsolution of (\ref{diri}) whereas the lower semicontinuous envelope of $W$ is a viscosity supersolution of (\ref{diri}). Then, the fact that $\underline{u}(\xi)= \psi(\xi)= \overline{u}(\xi)$ on $\partial\Omega$ and the comparison principle (Theorem \ref{cp} with $H\equiv 0$), guarantee that $W$ is continuous in $\overline\Omega$ and it solves  (\ref{diri}) in the viscosity sense.$\hfill\Box$
 
\subsection{Second case: more general operators with first order terms}
The next results provide the existence of a lower and an upper barrier for the following Dirichlet problem involving  operators depending on second and also on {\sl  first order} terms which degenerate at the characteristic points of the boundary. \\
Let us first consider the model  spherically symmetric example which will be useful to treat the general case providing appropriate local barrier functions for more general domains.\\
Take $0<R_1\leq R_2<\infty$ and  consider (\ref{0}) in an annular domain:
\begin{equation}\label{diri2}
\left\{
\begin{array}{ll}
F(\xi, D^2_{H^n}u)+H(\xi,\nabla_{H^n}u)=0\qquad &\mbox{in }\Omega= B^H_{R_2}(\eta_0)\setminus B^H_{R_1}(\eta_0)\\
u=0&\mbox{on }\partial (B^H_{R_2}(\eta_0)\setminus B^H_{R_1}(\eta_0))\,,
\end{array}
\right .
\end{equation}
 where  $F$ verifies (\ref{sigmaue}) and (\ref{sigmaf00}), $H$ is a continuous function  satisfying:
\begin{equation}\label{H}
|H(\xi,\nabla_{H^n}u)|\leq K|\nabla_{H^n}\rho(\eta_0^{-1}\circ\xi)||\nabla_{H^n}u| + M|\nabla_{H^n}\rho(\eta_0^{-1}\circ\xi)|^2\,.
\end{equation}
Since the characteristic points of $\partial (B^H_{R_2}(\eta_0)\setminus B^H_{R_1}(\eta_0))$ are $(\eta_{0_1},\dots,\eta_{0_{2n}},\eta_{0_{2n+1}}\pm {R^2_{i}})$, for $i=1\,,2$ and at these points the $\nabla_{H^n}\rho(\eta_0^{-1}\circ \xi)=0$, we deduce that the first order term degenerates at these points. \\
Let us point out that in the general setting  (\ref{sigmaf00}) and  (\ref{H}) $u=0$ is not a barrier (neither a lower nor an upper barrier).

\begin{lem}\label{barrierdiri2}
Let $F$ satisfy (\ref{sigmaue}) and (\ref{sigmaf00}) and $H$ satisfy  (\ref{H}). Then, there exist a lower barrier $\underline u$ and an upper barrier $\overline{u}$ for the problem (\ref{diri2}) which satisfy $\underline{u}(\xi)=\overline{u}(\xi)=0$ on $\partial\Omega$.
 \end{lem}
\noindent {\bf Proof}\quad
In order to construct  an  upper barrier for (\ref{diri2}) it is enough to choose a viscosity  supersolution $\overline u_2$  of $F+H\geq 0$ in $B^H_{R_2}(\eta_0)$ such that $\overline u_2>0$ in $B^H_{R_2}(\eta_0)$ and $\overline u_2=0$ on $\partial B^H_{R_2}(\eta_0)$, to take a  supersolution $\overline u_1$  of $F+H\geq 0$ in  $\left(B^H_{R_1}(\eta_0)\right)^c$ such that $\overline u_1>0$ in $\left(B^H_{R_1}(\eta_0)\right)^c$ and $\overline u_1=0$ on $\partial B^H_{R_1}(\eta_0)$,  and then to set $\overline u=\inf\{\overline u_1\,,\,\overline u_2\}$.

 Let us first construct   $\overline u_2$. The assumptions (\ref{sigmaue}),  (\ref{sigmaf00}) and  (\ref{H}) allow to look at the supersolution for  the problem:
\begin{equation}\label{supdiri2}
\left\{
\begin{array}{ll}
\tilde{\cal M}^-_{\lambda ,\Lambda}( D^2\overline{u}_2)-K|\nabla_{H^n}\rho(\eta_0^{-1}\circ\xi)|||\nabla_{H^n}\overline{u}_2|-M|\nabla_{H^n}\rho(\eta_0^{-1}\circ\xi)|^2\geq 0\quad &\mbox{in }B^H_{R_2}(\eta_0)\,,\\
\overline{u}_2>0 &\mbox{in } B^H_{R_2}(\eta_0)\,,\\
\overline{u}_2=0&\mbox{on }\partial B^H_{R_2}(\eta_0)\,.
\end{array}
\right .
\end{equation}
At this purpose, let us consider the following function:
\begin{equation}\label{ubar2}
\overline{u}_2(\xi)=\beta( e^{\alpha\frac{R_2^2}{2}}-e^{\alpha\frac{\rho^2(\eta_0^{-1}\circ \xi)}{2}})\,,
\end{equation}
where  $\beta\,,\alpha >0$ will be chosen later.

 Then $\overline{u}_2=0$ on $\partial B^H_R(\eta_0)$ and $\overline{u}_2>0$ in $B^H_R(\eta_0)$. Moreover, from Lemma \ref{lemmaeigen}, and the invariances of the vector fields $X_i$ with respect to the action $\circ$, it follows that the eigenvalues of $D^2_{H^n}\overline{u}_2$ are:
$$\begin{array}{l}
-\beta\alpha (1+\alpha\rho^2(\eta_0^{-1}\circ\xi))e^{\alpha\frac{\rho^2(\eta_0^{-1}\circ\xi)}{2}}|\nabla_{H^n}\rho(\eta_0^{-1}\circ\xi)|^2\mbox{ which is simple\,;}\\
-3\beta\alpha |\nabla_{H^n}\rho(\eta_0^{-1}\circ\xi)|^2 e^{\alpha\frac{\rho^2(\eta_0^{-1}\circ\xi)}{2}}\mbox{ which is simple\,;}\\
-\beta\alpha |\nabla_{H^n}\rho(\eta_0^{-1}\circ\xi)|^2 e^{\alpha\frac{\rho^2(\eta_0^{-1}\circ\xi)}{2}}\mbox{ with multiplicity $2n-2$\,.}
\end{array}
$$
Thus, 
$$
\begin{array}{ll}
\tilde{\cal M}^-_{\lambda ,\Lambda}(D^2\overline{u}_2)-K&|\nabla_{H^n}\rho(\eta_0^{-1}\circ\xi)||\nabla_{H^n}\overline{u}_2|-M|\nabla_{H^n}\rho(\eta_0^{-1}\circ\xi)|^2= \\
|\nabla_{H^n}\rho(\eta_0^{-1}\circ\xi)|^2 &\left [ \lambda \left(\beta\alpha (1+\alpha\rho^2(\eta_0^{-1}\circ\xi))e^{\alpha\frac{\rho^2(\eta_0^{-1}\circ\xi)}{2}}+ (2n+1)\beta\alpha e^{\alpha\frac{\rho^2(\eta_0^{-1}\circ\xi)}{2}}\right)\right.\\
&\left. -K\beta\alpha\rho(\eta_0^{-1}\circ\xi)e^{\alpha\frac{\rho^2(\eta_0^{-1}\circ\xi)}{2}} -M \right ]\,.
\end{array}
$$
Hence, being $e^{\alpha\frac{\rho^2(\eta_0^{-1}\circ\xi)}{2}}\geq 1$, in order to satisfy (\ref{supdiri2}) it is enough to choose  $\alpha\,,\,\beta>0$ in such a way that:
$$\lambda\beta\alpha (Q+\alpha \rho^2(\eta_0^{-1}\circ\xi))-K\beta\alpha\rho(\eta_0^{-1}\circ\xi)\geq M\,.$$
Then, taking $\alpha $  be such that 
$$\inf_{\rho}[\lambda Q+\lambda\alpha\rho^2-K\rho ]=C>0\,,$$
that is $\alpha=\frac{K^2}{2\lambda^2Q}$, and  $\beta\geq \frac{M}{\alpha C}$, (\ref{supdiri2}) holds true.

In order to construct $\overline{u}_1$ it is enough to take the supersolution of 
\begin{equation}\label{supdiri3}
\left\{
\begin{array}{ll}
\tilde{\cal M}^-_{\lambda ,\Lambda}( D^2\overline{u}_1)-K|\nabla_{H^n}\rho(\eta_0^{-1}\circ\xi)|||\nabla_{H^n}\overline{u}_1|-M|\nabla_{H^n}\rho(\eta_0^{-1}\circ\xi)|^2\geq 0\quad &\mbox{in }(B^H_{R_1}(\eta_0))^c\,,\\
\overline{u}_1>0 &\mbox{in } (B^H_{R_1}(\eta_0))^c\,,\\
\overline{u}_1=0&\mbox{on }\partial B^H_{R_1}(\eta_0)\,.
\end{array}
\right .
\end{equation}
At this purpose, let us take
\begin{equation}\label{ubar1}
\overline{u}_1(\xi)=\beta( e^{\frac{\alpha}{2R_1^2}}-e^{\frac{\alpha}{2\rho^2(\eta_0^{-1}\circ \xi)}})\,.
\end{equation}
 Then $\overline{u}_1=0$ on $\partial B^H_{R_1}(\eta_0)$ and $\overline{u}_1>0$ in $\Omega$. Moreover the eigenvalues of $D^2_{H^n}\overline{u}_1$ are:
$$\begin{array}{l}
-\frac{\beta\alpha}{\rho^4(\eta_0^{-1}\circ\xi)} (3+\frac{\alpha}{\rho^2(\eta_0^{-1}\circ\xi)})|\nabla_{H^n}\rho(\eta_0^{-1}\circ\xi)|^2 e^{\frac{\alpha}{2\rho^2(\eta_0^{-1}\circ\xi)}}\mbox{ which is simple\,;}\\
3\frac{\beta\alpha}{\rho^4(\eta_0^{-1}\circ\xi)} |\nabla_{H^n}\rho(\eta_0^{-1}\circ\xi)|^2 e^{\frac{\alpha}{2\rho^2(\eta_0^{-1}\circ\xi)}}\mbox{ which is simple\,;}\\
\frac{\beta\alpha}{\rho^4(\eta_0^{-1}\circ\xi)}  |\nabla_{H^n}\rho(\eta_0^{-1}\circ\xi)|^2 e^{\frac{\alpha}{2\rho^2(\eta_0^{-1}\circ\xi)}}\mbox{ with multiplicity $2n-2$\,.}
\end{array}
$$
Hence,
$$
\begin{array}{ll}
\tilde{\cal M}^-_{\lambda ,\Lambda}( D^2\overline{u}_1)-K&|\nabla_{H^n}\rho(\eta_0^{-1}\circ\xi)||\nabla_{H^n}\overline{u}_1|-M|\nabla_{H^n}\rho(\eta_0^{-1}\circ\xi)|^2=\\
  |\nabla_{H^n}\rho(\eta_0^{-1}\circ\xi)|^2&\left\{-\Lambda \beta\alpha \frac{(2n+1)}{\rho^4(\eta_0^{-1}\circ\xi)}e^{\frac{\alpha}{2\rho^2(\eta_0^{-1}\circ\xi)}}+\lambda\frac{\beta\alpha}{\rho^4(\eta_0^{-1}\circ\xi)} (3+\frac{\alpha}{\rho^2(\eta_0^{-1}\circ\xi)})e^{\frac{\alpha}{2\rho^2(\eta_0^{-1}\circ\xi)}}\right.\\
& \left. -K\frac{\beta\alpha}{\rho^3(\eta_0^{-1}\circ\xi)}e^{\frac{\alpha}{2\rho^2(\eta_0^{-1}\circ\xi)}}-M\right\}\geq 0
 \end{array}
 $$
 if 
\begin{equation}\label{condu1}
\frac{e^{\frac{\alpha}{2\rho^2(\eta_0^{-1}\circ\xi)}}}{\rho^4(\eta_0^{-1}\circ\xi)}\beta\alpha [-\Lambda (2n+1)+\lambda (3+\frac{\alpha}{\rho^2(\eta_0^{-1}\circ\xi)})-K\rho(\eta_0^{-1}\circ\xi)]\geq M
\end{equation}
and being $e^{\frac{\alpha}{2\rho^2(\eta_0^{-1}\circ\xi)}}\geq 1$, and $\rho(\eta_0^{-1}\circ\xi)\leq R_2$, (\ref{condu1}) holds true provided $\alpha$ satisfy $\beta\alpha\geq 1$ and
$$
\lambda (3+\frac{\alpha}{R_2^2})- KR_2-\Lambda(2n+1)\geq MR_2^4
$$
that is 
$$
\alpha\geq R_2^2(\frac{KR_2+\Lambda (2n+1)+MR_2^4}{\lambda}-3)\,.
$$
In order to provide a lower barrier, it is enough to take $\underline{u}=-\overline{u}$ which satisfies 
$$\left\{
\begin{array}{ll}
\tilde{\cal M}^+_{\lambda ,\Lambda}( D^2\underline{u})+K|\nabla_{H^n}\rho(\eta_0^{-1}\circ\xi)| | \nabla_{H^n}\underline{u}|+M|\nabla_{H^n}\rho(\eta_0^{-1}\circ\xi)|^2 \leq 0\qquad &\mbox{in }\Omega\,,\\
\underline{u}<0 &\mbox{in } \Omega\,,\\
\underline{u}=0&\mbox{on }\partial \Omega\,.
\end{array}
\right.
$$
The claim is accomplished. $\hfill\Box$
 
Analogously to the existence theorem for problem (\ref{diri}), by using Lemma \ref{barrierdiri2}, and the comparison principle (Theorem \ref{cp}), it results:
\begin{thm}\label{existencediri2} Let $F$ and $H$ satisfy the hypotheses of Lemma \ref {barrierdiri2}, then there exists a unique viscosity solution of problem (\ref{diri2}) which is continuous up to the boundary of $B^H_{R_2}(\eta_0)\setminus B^H_{R_1}(\eta_0)$.
\end{thm}
As it is immediate to verify by following the proof of Lemma \ref{barrierdiri2}, the results holds true also in the case where $R_1=0$. Indeed, in this case, the function $\overline u_2$ defined in (\ref{ubar2}) is an upper barrier and $- \overline u_2$ is a lower barrier. Hence we get the following existence theorem for the Heisenberg-ball:
\begin{thm}\label{existencediri3} Let $F$ and $H$ satisfy the hypotheses of Lemma \ref {barrierdiri2}, then there exists a unique viscosity solution to the problem
\begin{equation}\label{diri7}
\left\{
\begin{array}{ll}
F(\xi, D^2_{H^n}u)+H(\xi,\nabla_{H^n}u)=0\qquad &\mbox{in }B^H_{R}(\eta_0)\\
u=0&\mbox{on }\partial (B^H_{R}(\eta_0)\,,
\end{array}
\right .
\end{equation}
which is continuous up to the boundary of $B^H_{R}(\eta_0)$.
\end{thm}
Now let us consider the case of an open bounded set $\Omega$ satisfying (\ref{regolarephi2}) and the following condition:
\begin{equation}\label{regolarephi1}
\Omega \mbox{ verifies (\ref{sferapt}) at all characteristic points of  $\partial\Omega$}\,.
\end{equation}
Under the previous hypotheses  for $\Omega$, we wish to investigate the existence and uniqueness of viscosity solutions for the problem 
\begin{equation}\label{diri4}
\left\{
\begin{array}{ll}
F(\xi, D^2_{H^n}u)+H(\xi,\nabla_{H^n}u)=0\qquad &\mbox{in }\Omega\\
u=\psi(\xi)&\mbox{on }\partial \Omega\,,
\end{array}
\right .
\end{equation}
 continuous on ${\overline{\Omega}}$.

To prove the existence of viscosity solutions for the problem (\ref{diri4}), we need some hypotheses on the first order term $H$. Precisely, we assume the following conditions:
\begin{equation}\label{1condH}
|H(\xi,\sigma(\xi)p)|\leq K|\sigma(\xi)p|+M \qquad  \mbox{for all }\xi\in{\overline{\Omega}}\,,p\in\R^{2n+1}
\end{equation}
and for every characteristic point $\xi_0\in\partial\Omega$, there exists  $R>0$ such that:
\begin{equation}\label{hgen}
\begin{array}{l}
|H(\xi,\sigma(\xi)p)|\leq K|\nabla_{H^n}\Phi(\xi)||\sigma(\xi)p|+M|\nabla_{H^n}\Phi(\xi)|^2\qquad \forall \xi\in B(\xi_0,R)\cap\Omega\,,\forall p\in\R^{2n+1}

\end{array}
\end{equation}
where $B(\xi_0,R)$ denotes the euclidean ball with centre at  $\xi_0$ and radius $R$.\\
We point out that, as remarked by Bardi \& Mannucci in remark 7.1
of \cite{bardiman}, our result
under condition (\ref{hgen})  is  somehow "sharp", indeed if the first
order term has a "wrong" sign at the characteristic points,
they expect that the viscosity solution does not attain the boundary
datum (see references given in \cite{bardiman}).\\
Let us remind that the continuity of the viscosity solution to (\ref{diri4}) is a consequence of the existence of a lower and an upper local barrier at every point of the boundary (see \cite{bardiman,bardibottacin}). Precisely,
\begin{defi}Let $\xi_0\in\partial\Omega$. We say that $w$ is a lower (resp. upper) barrier at the point $\xi_0$ if it is a viscosity subsolution (resp. supersolution) of $F+H = 0$ in $\Omega$ and satisfy both $w\leq \psi$ (resp.  $w\geq \psi$ ) on $\partial\Omega$ and $\lim_{\xi\to\xi_0}w(\xi)=\psi(\xi_0)$.
\end{defi}
\begin{thm}\label{bardiman1} Under the hypotheses  (\ref{sigmaue}) and (\ref{sigmaf00}), (\ref{1condH})  and   (\ref{hgen}), if for all $\xi_0\in\partial\Omega$, there exist a local upper barrier and  a local lower barrier, then there exists a unique viscosity solution to (\ref{diri4}) which is continuous on $\overline\Omega$\,.
\end{thm}
For the proof see \cite{bardiman} and the references therein.\\
Therefore, in order to apply Theorem \ref{bardiman1}, we need to find a local upper and a local lower barrier at all points $\xi_0\in\partial\Omega$.
The strategy of the proof of the existence of local barriers at $\xi_0$  is different according to the case where $\xi_0$ is a  characteristic point of the boundary or not. 
In the second case the method follows essentially the ideas contained in \cite{bardiman} but for the sake of completeness we shall give all the details. 
On the contrary, in the first case we shall construct  local barriers using the ones found in the model case treated in Lemma \ref{barrierdiri2}. A key ingredient to adapt the barriers previously obtained for the spherically symmetric model case to this general one is the following  lemma:
\begin{lem}\label{characteristic}
Let $\Omega$ verify (\ref{regolarephi2}) and (\ref{regolarephi1}).
Then at  every characteristic point $\xi_0\in\partial\Omega$  (see (\ref{defcar})), it results:
\begin{equation}\label{compatib}
\lim_{\begin{array}{l}\xi\to\xi_0\\ \xi\in\overline{\Omega}\end{array}}\frac{|\nabla_{H^n}\Phi(\xi)|}{|\nabla_{H^n}\rho(\eta_0^{-1}\circ\xi)|}\leq C\,,
\end{equation}
where $\eta_0$ is the centre of the external Heisenberg--ball at $\xi_0$.
\end{lem}
\noindent {\bf Proof}\quad
Without loss of generality we can assume that the characteristic point is the origin. Morerover, if the external Heisenberg--ball at $0$ is not vertical (that is its centre $\eta_0$ does not lie on the $\xi_{2n+1}-$axis) then the condition (\ref{compatib}) immediately holds true since the origin is not characteristic for the Heisenberg--ball and then $|\nabla_{H^n}\rho(\eta_0^{-1}\circ\xi)|\ne 0$ in a neighborhood of the origin. 

Thus, let us assume that the centre $\eta_0=(0,\dots ,0\dots, t_0)$ belongs to the $\xi_{2n+1}-$axis.
Then being $\nabla \Phi(0)\not=0$ and $\nabla_{H^n} \Phi(0)=0$, by the Implicit Function Theorem, in a neighborhood of the origin, we have that $\Phi(\xi)=\xi_{2n+1}-h(\xi_{1},...,\xi_{n},\xi_{n+1},...,\xi_{2n})$, where $h$ has the same regularity as  $\Phi$.\\
Hence :

\[\begin{array}{ll}
\nabla_{H^n} \Phi(\xi)&=\nabla_{H^n} (\xi_{2n+1}-h(\xi_{1},...,\xi_{n},\xi_{n+1},...,\xi_{2n}))
\\
 &=\nabla_{H^n} (\xi_{2n+1})-\nabla_{H^n} (h(\xi_{1},...,\xi_{n},\xi_{n+1},...,\xi_{2n}))
\end{array}
\]
that means (for $i=1,...,n$)
\begin{equation}\label{heisf2}
\begin{array}{l}
X_i(\Phi(\xi))\!=\!2\xi_{i+n}-
\frac{\partial}{\partial \xi_i} h(\xi_{1},...,\xi_{n},\xi_{n+1},...,\xi_{2n}),\\
X_{i+n}(\Phi(\xi))\!=\!-2\xi_i
-\frac{\partial}{\partial \xi_{i+n}} h(\xi_{1},...,\xi_{n},\xi_{n+1},...,\xi_{2n}).
\end{array}
\end{equation}
Now, since $\nabla_{H^n} \Phi(0)=0$, we get for  $j=1,...,2n$
\[
\frac{\partial h}{\partial \xi_{j}}(0)=0.
\]
By using a Taylor expansion of $\nabla h$ at $0$ we obtain
\begin{equation}
\label{T1}
\vert \nabla_{H^n} \Phi(\xi)\vert =2(\sum_{i=1}^{2n}\xi_i^2)^{\frac 1 2}+O((\sum_{i=1}^{2n}\xi_i^2)^{\frac 1 2}),
\end{equation}
On the other hand
\begin{equation}
\label{T2}
\vert \nabla_{H^n} \rho(\eta_0^{-1}\circ \xi)\vert =\frac{(\sum_{i=1}^{2n}\xi_i^2)^{\frac 1 2}}{((\sum_{i=1}^{2n}\xi_i^2)^2+(\xi_{2n+1}-t_0)^2)^{\frac 1 4}}\,.
\end{equation}
Then,  (\ref{T1}) and (\ref{T2}) yield:
\begin{equation}
\label{T3}
\frac{\vert \nabla_{H^n} \Phi(\xi)\vert}{\vert \nabla_{H^n} \rho(\eta_0^{-1}\circ \xi)\vert }=\frac
{(2(\sum_{i=1}^{2n}\xi_i^2)^{\frac 1 2}+O((\sum_{i=1}^{2n}\xi_i^2)^{\frac 1 2}))((\sum_{i=1}^{2n}\xi_i^2)^2+(\xi_{2n+1}-t_0)^2)^{\frac 1 4}
}
{(\sum_{i=1}^{2n}\xi_i^2)^{\frac 1 2}}.
\end{equation}

Therefore $\frac{\vert \nabla_{H^n} \Phi(\xi)\vert}{\vert \nabla_{H^n} \rho(\eta_0^{-1}\circ \xi)\vert }$ remains bounded in a neighborhood of the origin and the claim follows. $\hfill\Box$

\begin{esempio} If $\Omega$ is a Heisenberg--ball, and we assume that the characteristic point is the origin. Then the external Heisenberg ball has its centre at $\eta_0=(0,\dots,0,\dots,t_0)$. Assuming $t_0<0$, we get that $\Omega=\{t>h(\xi_1,\dots,\xi_{2n})\}$ where 
$$h(\xi_1,\dots,\xi_{2n})=t_0\left(1-\sqrt{1-\frac{(\sum_{i=1}^{2n}\xi_i^2)^2}{t_0^2}}\right)\,.$$
In this case with an easy calculation we obtain that
$$
\lim_{\begin{array}{l}\xi\to 0\\ t\to 0\end{array}}\frac{\vert \nabla_{H^n} \Phi(\xi)\vert}{\vert \nabla_{H^n} \rho(\eta_0^{-1}\circ \xi)\vert }=\sqrt{|t_0|+\frac{1}{|t_0|}}\,.
$$
Let us point out that if $t_0$ tends to zero (that is the external Heisenberg--ball degenerates to a point)  the limit becomes infinity.$\hfill\Box$
\end{esempio} 
Now, let us construct  a local upper and a local lower barrier at all points $\xi_0\in\partial\Omega$. \\
This is the claim of the following lemma: 
\begin{lem}\label{barrierdiri4}
Let $\Omega$ be a bounded domain verifying (\ref{regolarephi2}) and (\ref{regolarephi1}). Let 
 $\psi$ be a continuous  function on $\partial\Omega$. 
Let $F$ satisfy (\ref{sigmaue}) and (\ref{sigmaf00}) and $H$ satisfy (\ref{1condH})  and (\ref{hgen}). Then, there exist a local lower barrier $\underline u$ and a local upper barrier $\overline{u}$ for the problem 
(\ref{diri4}) at all points $\xi_0\in\partial\Omega$.
 \end{lem}
\noindent {\bf Proof}\quad

\noindent {\it First step} \\
Let us  first construct a supersolution $w_1$ and a subsolution $w_2$ of the equation in (\ref{diri4}) which respectively satisfy $w_1\geq \psi$ and  $w_2\leq \psi$ on $\partial\Omega$. 

At this purpose, we can follow the same idea of Bardi\&Mannucci \cite{bardiman}. Indeed the function 
$$w_1(\xi)=k\left(\beta-\exp(\frac{\mu|\xi|^2}{2})\right)\qquad \mbox{ with }\beta\geq\sup_{\partial\Omega}\left(|\psi(\xi)|+\exp(\frac{\mu|\xi|^2}{2})\right)$$
is a viscosity supersolution of $F+H\geq 0$, for $k\geq 1$ and $\mu$ sufficiently large. Indeed, by the assumption (\ref{1condH})  we get that 
$$
F(\xi, D^2_{H^n}w_1)+H(\xi,\nabla_{H^n}w_1)\geq \tilde{\cal M}^-_{\lambda ,\Lambda}( D^2{w}_1)-K|\nabla_{H^n}{w}_1|-M \geq 0\,,
$$
as it is easy to check by taking into account that:
$$\begin{array}{l}
\nabla_{H^n}w_1=-k\mu\exp(\frac{|\xi|^2}{2})\nabla_{H^n}\left(\frac{|\xi|^2}{2}\right)\\
D^2_{H^n}w_1=-k\exp(\mu\frac{|\xi|^2}{2})[\mu^2\nabla_{H^n}\left(\frac{|\xi|^2}{2}\right)\otimes\nabla_{H^n}\left(\frac{|\xi|^2}{2}\right)+\sigma(\xi)\sigma^T(\xi)]\,,
\end{array}
$$
which implies that
$$\begin{array}{l}
\tilde{\cal M}^-_{\lambda ,\Lambda}( D^2w_1)\geq k\exp(\frac{\mu|\xi|^2}{2})\left (\mu^2{\cal M}^-_{\lambda ,\Lambda}(-\nabla_{H^n}\left(\frac{|\xi|^2}{2}\right)\otimes\nabla_{H^n}\left(\frac{|\xi|^2}{2}\right)) +{\cal M}^-_{\lambda ,\Lambda}(-\sigma(\xi)\sigma^T(\xi))\right)\\
=k\lambda \exp(\frac{|\xi|^2}{2})\left(\mu^2|\nabla_{H^n}\left(\frac{|\xi|^2}{2}\right)|^2+2n+4(\sum_{i=1}^{2n}\xi^2)\right)\,.
\end{array}
$$
Hence,
$$
\begin{array}{l}
F(\xi, D^2_{H^n}w_1)+H(\xi,\nabla_{H^n}w_1)\geq \\
k \exp(\frac{|\xi|^2}{2})[\mu^2\lambda|\nabla_{H^n}\left(\frac{|\xi|^2}{2}\right)|^2-K\mu|\nabla_{H^n}\left(\frac{|\xi|^2}{2}\right)|+\lambda(4(\sum_{i=1}^{2n}\xi^2)+2n)]-M\geq 0
\end{array}
$$
if $\mu$ is sufficiently large.

 Then, taking $w_2(\xi)=-w_1(\xi)$ we obtain a subsolution.

\noindent {\it Second step}\\
Now let us construct a local upper barrier in a neighborhood of a  point $\xi_0\in\partial\Omega$. 

We need to follow different strategies according to  the fact that $\xi_0$ is a characteristic point or not.

Let us first consider the case $\xi_0\in\partial\Omega$ is not characteristic (that is $|\nabla_{H^n}\Phi(\xi_0)|\neq 0$). 

Then, by continuity we find a neighborhood of $\xi_0$ where   $|\nabla_{H^n}\Phi(\xi)|\neq 0$. Let us call $B(\xi_0;R)\cap\Omega$, this neighborhood. By using the same idea of Bardi\&Mannucci, \cite{bardiman}, consider the function 
$$W(\xi)=1-\exp{\left(-\mu(\Phi(\xi)+\frac{\alpha}{2}|\xi-\xi_0|^2)\right)}\,.$$
Obviously, $W(\xi_0)=0$, $W(\xi)\geq 0$ on $\partial\Omega$ and moreover, 
$$
\nabla_{H^n}W(\xi_0)=\mu\nabla_{H^n}\Phi(\xi_0)\,,
$$
$$D^2_{H^n}W(\xi_0)=\mu\left[D^2_{H^n}\Phi(\xi_0)-\mu(\nabla_{H^n}\Phi(\xi_0)\otimes\nabla_{H^n}\Phi(\xi_0))+\alpha(\sigma(\xi_0)\sigma^T(\xi_0))\right]\,.
$$
Hence,
$$\begin{array}{l}
\tilde{\cal M}^-_{\lambda ,\Lambda}( D^2W(\xi_0))\geq \mu\tilde{\cal M}^-_{\lambda ,\Lambda}( D^2\Phi(\xi_0))-\mu^2\tilde{\cal M}^-_{\lambda ,\Lambda}(\nabla_{H^n}\Phi(\xi_0)\otimes\nabla_{H^n}\Phi(\xi_0) )+\mu\alpha\tilde{\cal M}^-_{\lambda ,\Lambda}( \sigma(\xi_0)\sigma^T(\xi_0))\\
=\mu^2(-\tilde{\cal M}^-_{\lambda ,\Lambda}(\nabla_{H^n}\Phi(\xi_0)\otimes\nabla_{H^n}\Phi(\xi_0) )+\frac{1}{\mu}\tilde{\cal M}^-_{\lambda ,\Lambda}( D^2\Phi(\xi_0))+\frac{\alpha}{\mu}\tilde{\cal M}^-_{\lambda ,\Lambda}( \sigma(\xi_0)\sigma^T(\xi_0)))\\
\geq\frac{\mu^2}{2}(\Lambda|\nabla_{H^n}\Phi(\xi_0)|^2)>0 \quad\mbox{ if $\mu$ is large enough .}\end{array}
$$
Therefore, by choosing $\mu$ sufficiently large,
$$
\begin{array}{l}
\tilde{\cal M}^-_{\lambda ,\Lambda}( D^2W(\xi_0))-K|\nabla_{H^n}W(\xi_0))|-M\\
\geq (\frac{\mu^2}{4}\Lambda|\nabla_{H^n}\Phi(\xi_0))|^2>0
\end{array}
$$
and, by continuity,
$$
\begin{array}{l}
\tilde{\cal M}^-_{\lambda ,\Lambda}( D^2W(\xi ))-K|\nabla_{H^n}W(\xi))|-M\geq 0\qquad\forall \xi \in B(\xi_0;R)\cap\Omega\,.
\end{array}
$$

Now let $\xi_0\in\partial\Omega$ be a characteristic point (that is $|\nabla_{H^n}\Phi(\xi_0)|= 0$). 

Let $\eta_0\in\Omega^c$ be the centre of the Heisenberg--ball touching $\partial\Omega$ at $\xi_0$. Take the ball with centre at $\eta_0$ and radius $r_1=\rho (\eta_0^{-1}\circ\xi_0)$, where $\rho$ is defined in (\ref{sferaext}).
The hypotheses (\ref{hgen}) and (\ref{compatib}) yield 
$$\begin{array}{l}
F(\xi, D^2_{H^n}W)+H(\xi,\nabla_{H^n}W)\geq 
\tilde{\cal M}^-_{\lambda ,\Lambda}( D^2W)-
K|\nabla_{H^n}\Phi(\xi)||\nabla_{H^n}W|-M|\nabla_{H^n}\Phi(\xi)|^2\\
\geq 
\tilde{\cal M}^-_{\lambda ,\Lambda}( D^2W)-
{\tilde K}|\nabla_{H^n}\rho(\eta_0^{-1}\circ\xi)||\nabla_{H^n}W|-{\tilde M}|\nabla_{H^n}\rho(\eta_0^{-1}\circ\xi)|^2\geq 0\mbox{ in } \Omega\cap B(\xi_0,R)
\end{array}
$$
for a suitable $R>0$. By choosing 
$W(\xi)\equiv \overline{u}_1$ defined in (\ref{ubar1}) of Lemma \ref{barrierdiri2} associated with the constants ${\tilde K}$ and ${\tilde M}$ instead of $K$ and $M$, we obtain 
$$F(\xi, D^2_{H^n}W)+H(\xi,\nabla_{H^n}W)\geq 0 \qquad\mbox{ in }\Omega\cap B(\xi_0,R)$$
for some $R>0$.  Since  $W$ vanishes on $\partial  B^H_{r_1}(\eta_0)$, it results $W(\xi_0)=0$, $W(\xi)>0$ in $\Omega\cap B(\xi_0,R)$.

\noindent{\it Third step }\\
The function
\begin{equation}\label{incolla}
w_{\eps,\xi_0}(\xi)=\left\{\begin{array}{ll}
\min{\left\{\psi(\xi_0)+\eps +\tau W(\xi)\,;\,w_1(\xi)\right\}}\qquad &\xi\in\overline{B(\xi_0;R)\cap\Omega}\\
w_1(\xi)&\mbox{otherwise}
\end{array}
\right.
\end{equation}
with $\tau$ large enough to get that $w_{\eps,\xi_0}(\xi)=w_1(\xi)$ when $|\xi-\xi_0|=R$ is a local upper barrier at $\xi_0$.
Now repeating this argument at all points $\xi_0$, and setting
$$\overline{u}(\xi):=\inf\{w_{\eps,\xi_0}(\xi)\,\:\,\eps>0, \,\xi_0\in\partial\Omega\,\}\mbox{   for }\xi\in\overline\Omega\,,
$$
we get an upper barrier for (\ref{diri4}), which attains the boundary datum $\psi$. By following the same strategy as in Lemmas \ref{barrierdiri} and \ref{barrierdiri2}, we can construct a lower barrier, and we conclude.$\hfill\Box$

So, the following result holds true:
\begin{thm}\label{existencediri5} Let $F$, $H$ and  $\Omega$ satisfy the hypotheses of Lemma \ref{barrierdiri4}. Then there exists a unique viscosity solution to the problem (\ref{diri4}) continuous up to the boundary of $\Omega$.
\end{thm}

Let us remark that for $H=0$  Theorem \ref{existencediri5}  provides an existence result with different regularity conditions on $\partial\Omega$ with respect to the ones of Theorem \ref{existence1}. Indeed, in Theorem \ref{existence1} we asked  the existence of the external Heisenerg--ball at all points of $\partial\Omega$ but we did not require (\ref{regolarephi2}), whereas in Theorem \ref{existencediri5}  we need  (\ref{regolarephi1}) and  (\ref{regolarephi2}).
\bigskip

\section{Qualitative properties of viscosity supersolutions}\label{qua}
The just constructed barriers for Pucci--Heisenberg operators allow to find some qualitative properties of   Hadamard and Liouville type, generalizing the ones found in \cite{acfl,cdc2} for the uniformly elliptic fully non linear case.

The first result is related to the behaviour of the minima on intrinsic balls of viscosity  supersolutions to the Pucci--Heisenberg operators. Then, we apply it in order to find some non-linear Liouville results and a weak Harnack inequality for radial supersolutions  with respect to $\rho$ of the inequality  $F(\xi,D^2_{H^n}u)\geq 0$.

\begin{thm}\label{hadamdeg} {\em (Nonlinear Degenerate Hadamard Theorems)}
Let $\Omega$ be a domain of $\R^{2n+1}$ containing the closed intrinsic ball ${\overline B}^H_R(0)$
centered at the origin and with radius $R >0$. Then:
\begin{itemize}
\item[(i)]  
 if $u\in LSC(\Omega)$ is a viscosity solution of
\[
\tilde{\cal M}^-_{\lambda ,\Lambda}(D^2 u)\geq 0\, \qquad\hbox{ in }\Omega\,,
\]
 then the function
\[
m(r)\, =\, \min_{\rho(\xi)\leq r}u(\xi)\,,\qquad r<R
\]
is, respectively, a concave function of $\log r$ if $\alpha =2$ and of $r^{2-\alpha }$ if $\alpha
\neq 2$, with $\alpha$ given by (\ref{sigmanewN1}).
 More precisely, for every fixed $r_1<R$,
it results
\begin{equation}\label{sigmaHad+-}
m(r)\geq \frac{G(r)}{G(R)}
m(R)+\left(1-\frac{G(r)}{G(R)}\right ) m(r_1) \;\;,\;\; \forall 
r\in [r_1, R]
\end{equation}
with
\begin{equation}\label{GsigmaHad+-}
G(r)=
\left\{
\begin{array}{ll}
\log (r/r_1)& \qquad
\mbox{if} \  \alpha =2 \nl
r^{2-\alpha }-r_1^{2-\alpha}& \qquad
\mbox{if} \  \alpha \neq 2\, . \\
\end{array}
\right.
\end{equation}
\item[(ii)]  
if $u\in LSC(\Omega)$ is a viscosity solution of
\[
\tilde{\cal M}^+_{\lambda ,\Lambda}(D^2 u)\geq 0\,  \qquad\hbox{ in }\Omega\,,
\]
then 
$m(r)$ is a concave function of  $r^{2-\beta }$, with $\beta$ given by (\ref{sigmanewN2}).
 More precisely, for every fixed $r_1<R$ it satisfies (\ref{sigmaHad+-}) with 
 \begin{equation}\label{GsigmaHad--}
G(r)=r^{2-\beta }-r_1^{2-\beta }\,.
\end{equation}
\end{itemize}
\end{thm}
Before giving the proof of the theorem, let us observe that, by the
relationship (\ref{sigmaeqmpm}) between  $\tilde{\cal M}^+_{\lambda
,\Lambda}$ and $\tilde{\cal M}^-_{\lambda ,\Lambda}$, statement  $(i)$ is
equivalent to the following one:
\begin{itemize} 
\item[(j)]
if $u\in USC(\Omega)$ is a viscosity solution of
\[
\tilde{\cal M}^+_{\lambda ,\Lambda}(D^2 u)\leq 0\,\qquad\hbox{ in }\Omega\,,
\]
 then the function $M(r)\, =\, \max_{\rho(\xi)\leq r}u(\xi)$ satisfies
\begin{equation}\label{sigmaHad-+}
M(r)\leq \frac{G(r)}{G(R)}
M(R)+\left(1-\frac{G(r)}{G(R)}\right ) M(r_1) \;\;,\;\; \forall 
r\in [r_1, R]
\end{equation}
with $G$ given by (\ref{GsigmaHad+-}).
\end{itemize}
Analogously, an equivalent form of $ (ii)$ is
\begin{itemize}
\item[(jj)]
if $u\in USC(\Omega)$ is a viscosity solution of
\[
\tilde{\cal M}^-_{\lambda ,\Lambda}(D^2 u)\leq 0\, \qquad\hbox{ in }\Omega\,,
\]
 then $M(r)$ satisfies (\ref{sigmaHad-+}) with $G$ given by (\ref{GsigmaHad--}).
\end{itemize}
\bigskip
\noindent {\bf Proof}\quad  By the assumptions, the respectively increasing
and decreasing functions $M(r)$ and $m(r)$ are well defined in $[0,R]$.

 Let us  consider the case $(i)$, that is, let $u\in LSC(\Omega)$ be a viscosity
solution of $\tilde{\cal M}^-_{\lambda ,\Lambda}(D^2 u)\geq 0$.\\
\noindent Fixed $0<r_1<R$, let $\psi (r)=-\varphi_1(r)$, with $\varphi_1(r)$  defined by
(\ref{sigmasolfon1}), with constants $C_1 \geq 0$ and $C_2 \in \R$ chosen in such a way
that $\psi (R)=m(R)$ and $\psi (r_1)=m(r_1)$. This yields:
\[
\psi (r)=\frac{G(r)}{G(R)}
m(R)+\left(1-\frac{G(r)}{G(R)}\right ) m(r_1)
\]
with $G$ given by  (\ref{GsigmaHad+-}).
We know that the function $\Psi (\xi)=\psi (\rho(\xi))$ is a viscosity solution of
equation (\ref{sigmaomo-}). Applying the Comparison Principle (Theorem \ref{cp} for $H=0$) to the functions $u(\xi)$ and 
$\Psi (\xi)$ in the intrinsic annular domain
$\{ r_1<\rho(\xi)<R\} \subset \Omega$, we deduce that \[
u(\xi)\geq \Psi (\xi) \qquad  \mbox{  in }\  \{ r_1\leq \rho(\xi)\leq R\} \, .
\]
Hence, $m(r)\geq \psi (r)$ for all $r$ in $[r_1,R]$ and the claim is proved.\\
The proof of $(ii)$ is completely analogous to that of $(i)$, with the
obvious difference that now $u$ has to be compared with the function
$\Phi_2(\xi)=\varphi_2 (\rho(\xi))$, where $\varphi_2$ is given by
(\ref{sigmasolfon2}). $\hfill\Box$

\bigskip 
   
\noindent Looking at the previous result, as well as at the
just constructed \lq\lq fundamental solutions" of  equations (\ref{sigmaomo+})
and (\ref{sigmaomo-}), we expect a
Liouville type theorem in the case  $\Omega=\R^{2n+1}$, when the \lq\lq fundamental solutions" go to infinity as $\rho\rightarrow +\infty$, that is  for
  bounded from below (above) viscosity supersolutions
(subsolutions) of (\ref{sigmaomo-}) ((\ref{sigmaomo+})) in all of $\R^{2n+1}$, when
the parameter $\alpha$, defined by (\ref{sigmanewN1}), satisfies $\alpha \leq 2$.
Indeed, in the other cases, by using the \lq\lq fundamental solutions" we can construct a non costant viscosity solution. As an example, in the case  $(ii)$, a strictly positive viscosity supersolution of 
\begin{equation}\label{eqcontroes}\tilde{\cal M}^+_{\lambda ,\Lambda}(D^2 u)\geq 0\, \qquad\mbox{ in }\R^{2n+1}
\end{equation}
is the function
\begin{equation}\label{controes}
 u=\left\{ \begin{array}{ll}
\rho^{2-\beta}(\xi)  &\mbox{if $\rho(\xi)\geq R$} \\
R^{2-\beta} & \mbox{if $\rho(\xi)<R$}\,,
\end{array}\right.
\end{equation}
where $R>0$ is arbitrarily fixed. Indeed, $\rho^{2-\beta}(\xi) $ is a classical solution of 
$$
\tilde{\cal M}^+_{\lambda ,\Lambda} (D^2 u)=0\mbox{   in }\R^{2n+1}\setminus\{0\}
$$
(see Lemma \ref{funsolm+}), whereas $R^{2-\beta}$ is a solution of the same equation in the whole 
space. Therefore, well known stability properties of viscosity solutions 
imply (see \cite{user}) that $u\equiv\min{\{R^{2-\beta};\rho^{2-\beta}(\xi)\}}$ is a viscosity solution of (\ref{eqcontroes}).\\
The following result holds true: 
 \begin{thm}\label{1} Let $u\in C(\R^{2n+1})$ be a viscosity solution either of
\begin{equation}\label{sigmam-}
 \tilde{\cal M}^-_{\lambda,\Lambda}(D^2u)\geq 0 \qquad\mbox{ in }\R^{2n+1}\, 
 \end{equation}
  or of
 \begin{equation}\label{sigmam+}
\tilde{\cal M}^+_{\lambda,\Lambda}(D^2u)\leq 0  \qquad\mbox{ in }\R^{2n+1}.
\end{equation}

 If $u$ is, respectively, bounded either from below or from above,  and if
the parameter $\alpha$, defined by (\ref{sigmanewN1}), satisfies $\alpha \leq 2$
\mbox{(}i.e.   $Q\leq \frac{\Lambda}{\lambda}+1$\mbox{)}, then $u$ is
constant.  \end{thm}
\noindent {\bf Proof} \quad Consider the case $ \tilde{\cal
M}^-_{\lambda,\Lambda}(D^2u)\geq 0$, $u$ bounded from below. 
By the previous theorem (case $(i)$), $u$ satisfies (\ref{sigmaHad+-}), with $G$ given by (\ref{GsigmaHad+-}) for every fixed $0<r_1 <R$.

\noindent Being $m(r)$ a bounded function since $u$ is bounded from 
below, and
being $\alpha \leq 2$, passing to the limit as $R\to +\infty$ in 
(\ref{sigmaHad+-})
leads to
\[
m(r)\geq m(r_1) \ \ \ \ \forall \, r\geq r_1 >0\, .
\]
Since $m(r)$ is obviously a decreasing function, we deduce
that $m(r)\equiv const.=m(0)=u(0)$. Therefore, $u$ attains its minimum at an
interior point and, by the Strong Maximum Principle stated in \cite{bdl2}, $u$ is constant.$\hfill\Box$

Another consequence of the Hadamard theorem is the following version 
of the weak Harnack inequality for radial supersolutions.

\noindent \begin{thm}\label{harnack} Let $u\in LSC( B^H_{2R}(0))$ be a radial viscosity solution of
\[
u\geq 0, \qquad F(\xi,D^2_{H^n}u)\geq 0\, \qquad\mbox{
in } B^H_{2R}(0)\, , 
\] 
with $F$ verifying (\ref{sigmaue}) and (\ref{sigmaf00}).
Then 
$u$ satisfies the following  weak Harnack inequality:
  \begin{equation}\label{har}
{\rm meas}\left( B^H_{\frac{R}{2}}\cap \{ u(\rho)>t\} \right) \leq
 \frac{C\, R^Q}{t^{\frac{Q}{\beta-2}}}\, \left( u(R)\right) ^{\frac{Q}{\beta-2}}\qquad \forall t>0\,,
\end{equation}
with $\{ u(\rho)>t\}=\{ \xi\in B^H_{2R}(0) \, :\,  u(\rho(\xi))>t\}$ and ${\rm meas}(E)$ equals to the
Lebesgue measure of the measurable set $E\subset \R^{2n+1}$.
\end{thm}

\noindent {\it Proof.} By (\ref{sigmastimafm-}), $u$ satisfies in the viscosity sense:
\[
u\geq 0, \qquad {\cal M}^+_{\lambda,\Lambda}(D^2u)\geq 0\, \qquad\mbox{
in } B^H_{2R}(0)\, .
\] 
Then, the strong maximum principle yields that the minimum
of the radial function $\min_{\xi\in {\overline B}^H_{R}}u(\rho (\xi))$ is attained at the boundary of the
set $B^H_{R}$. Moreover, by the Hadamard theorem we get that 
\begin{equation}\label{cresce} 
    m(\rho)\rho^{\beta -2}\mbox{ is an increasing function. }
    \end{equation}
    Hence, being in this case $m(\rho)=u(\rho)$, we get
  \begin{equation}\label{harmx}
\begin{array}{ll}\dis
{\rm meas}\left( B^H_{\frac{R}{2}}\cap \{m(\rho)>t\} \right) &\leq 
{\rm meas}\left( B^H_{\frac{R}{2}}\cap
\left\{\frac{m(R)R^{\beta-2}}{\rho^{\beta-2}}>t\right\} \right)\\
\dis &\leq
 \frac{C\, R^Q}{t^{\frac{Q}{\beta-2}}}\, \left( m(R)\right) ^{\frac{Q}{\beta-2}}\qquad \forall t>0\,.
\end{array}
\end{equation}
This concludes the proof.
\hfill$\Box$
\par\medskip


 
\end{document}